\documentclass[onefignum,onetabnum,hidelinks]{siamart190516}

\usepackage{lipsum}
\usepackage{amsfonts}
\usepackage{graphicx}
\usepackage{epstopdf}
\usepackage{algorithmic}
\usepackage{color}
\ifpdf
  \DeclareGraphicsExtensions{.eps,.pdf,.png,.jpg}
\else
  \DeclareGraphicsExtensions{.eps}
\fi

\newsiamremark{remark}{Remark}
\newsiamremark{hypothesis}{Hypothesis}
\crefname{hypothesis}{Hypothesis}{Hypotheses}
\newsiamthm{claim}{Claim}

\title{Optimal Control of Parameterized Maxwell's System: Reduced Basis, Convergence Analysis, and A Posteriori Error Estimates\thanks{This work is partially supported by the Air Force Office of Scientific Research (AFOSR) under Award NO: FA9550-19-1-0036 and NSF grants DMS-1818772 
and DMS-1913004.}
}

\author{
Harbir Antil\thanks{The Center of Mathematics and Artificial Intelligence (CMAI) and Department of
Mathematical Sciences, George Mason University, Fairfax, VA 22030, USA
(\email{hantil@gmu.edu}).}
\and
Tran Nhan Tam Quyen\thanks{Institute for Numerical and Applied Mathematics, University of Goettingen, Lotzestr. 16-18, 37083 Goettingen, Germany
(\email{quyen.tran@uni-goettingen.de}).}
}

\usepackage{amsopn}

\ifpdf
\hypersetup{
}
\fi

\def\mb{\mathbf}

\usepackage[textsize=small]{todonotes}
\begin{document}

\maketitle

% REQUIRED
\begin{abstract} 
We consider control constrained optimal control problems governed by parameterized stationary Maxwell's system with the Gauss's law. The parameters enter through dielectric, magnetic permeability, and charge density. Moreover, the parameter set is assumed to be compact. We discretize the electric field by a finite element method and use variational discretization concept to discretize the control. We create a reduced basis method for the optimal control problem and establish uniform convergence of the reduced order solutions to that of the original high dimensional problem provided that the snapshot parameter sample is dense in the parameter set, with an appropriate parameter separability rule. Finally, we establish the absolute a posteriori error estimator for the reduced order solutions and the corresponding cost functions in terms of the state and adjoint residuals. 

\end{abstract}

%In this paper we investigate the reduced order solution of the optimal control problem governed by a parameterized stationary Maxwell system with the Gauss law. In this context the dielectric, the magnetic permeability and the charge density are assumed to be known, where the control is constrained of general type and the parameter set is compact. We approximate the electric field of the Maxwell system in finite element spaces. Adopting the variational discretization concept, we consider a weighted parameterized optimal control problem. Utilizing techniques from the primal reduced basis approach, we construct a reduced basis surrogate model for the aforementioned optimal control problem. We prove the uniform convergence of reduced order solutions to that of the original high dimensional problem provided the snapshot parameter sample being dense in the parameter set and with an appropriate parameter separability rule. Furthermore, we establish the absolute a posteriori error estimator for the reduced order solutions and the corresponding cost functional which deals with the norm of the residuals of the state and adjoint equations.

% REQUIRED

% REQUIRED

%%%%%%%%%%%%%%%%%%%%%%%%%%%%%%%%%%%%%%%
\section{Introduction}
%%%%%%%%%%%%%%%%%%%%%%%%%%%%%%%%%%%%%%%

Maxwell's equations with the Gauss's law play a central role in many day-to-day applications. 
However, the underlying 
coefficients in these equations, such as dielectric, magnetic permeability, and charge density contains 
parameters which must be inferred from experiments or treated as random variables. In many cases, 
these parameterized equations must be queried  for different parameters, many times over and thus 
the problem quickly becomes intractable. This issue 
is only exasperated when dealing with optimization problems with such parameterized equations as 
constraints. The goal of this paper is to create efficient numerical methods, using reduced basis 
method, to solve the optimization problems governed by stationary Maxwell's system with the 
Gauss's law as constraints. 

We discretize  these equations using a finite element method and carry out a variational discretization 
for the control. The finite element system for the PDE is a parameterized constrained saddle point system. 
It can be very expensive to solve, especially on fine meshes and for many parameter queries (cf.\  \cite{Hesthaven16,AQ16,Ha17,HaOh08}). From a reduced basis point of view, one needs a surrogate 
model for the system. Furthermore, since the reduced basis approach considers a suboptimal problem, 
convergence analysis and error estimates for the reduced order solution to that of the original 
high dimensional problem are crucial, which are investigated in the present paper.

For completeness, we mention that the optimal control problems governed by the \emph{non-parameterized} Maxwell systems have attracted a great deal of attention from many scientists in the last decades. For surveys on the subject, we refer the reader to, e.g., \cite{BoYo16,KoLa12,Ni14,yousept-zou} and the references therein. The construction of the reduced basis methods for parameterized Maxwell systems can be found in \cite{Be17,Che11,Che10,Ha15,HeBe13,Hes15}. Moreover, an incomplete list of references that considers the analysis of parameterized optimal control problems (not Maxwell) can be found in \cite{Ba16,de10,ItRa98,KaGr14,Ka18,Ne13,To11}. 

We in \S\ref{fem} first present numerical analysis for Maxwell systems with N\'ed\'elec finite elements. To the best of our knowledge, the optimal control of such a system is not known. The main results of our paper are contained in Theorem \ref{8-5-19ct17}, where we prove the uniform convergence of reduced order solution, to the optimal control problem, to that of the original high dimensional problem, and in Theorem \ref{2-5-19-DL} where we establish the absolute a posteriori error estimator for the reduced order solutions. Numerical implementation will be part of a future work. 

The remainder of the paper is organized as follows. In \Cref{s:prob_f}, we state the problem under consideration. 
\Cref{12-6-19ct2} is devoted to some functional spaces and the finite element method for the system \cref{9-2-18ct1}. Primal reduced basis approach for the optimal control problem and first order optimality conditions are presented in \Cref{RBM}. Convergence analysis and a posteriori error estimates for the reduced basis approximations are respectively discussed in \Cref{convergence} and \Cref{Post-estimate}.

%%%%%%%%%%%%%%%%%%%%%%%%%%%%%%%%%%%%%%%
\section{Problem Formulation}
\label{s:prob_f}
%%%%%%%%%%%%%%%%%%%%%%%%%%%%%%%%%%%%%%%

Let $\Omega$ be an open, bounded and connected set in $\mathbb{R}^3$ with the Lipschitz boundary $\partial \Omega$, and $\mathcal{P} \subset \mathbb{R}^p$ is a compact set of parameters. In this paper we deal with the following $\mu$-parameterized stationary Maxwell's system fulfilled by the electric field $\mathbf{E}$:
\begin{equation}\label{9-2-18ct1}
\begin{cases}
\nabla\times \left(\sigma^{-1}(\mathbf{x};\mu)\nabla\times \mathbf{E}(\mathbf{x};\mu)\right)&=~ \epsilon(\mathbf{x};\mu) \mathbf{u}(\mathbf{x}), \quad(\mathbf{x};\mu) \in \Omega\times \mathcal{P},\\
\nabla\cdot(\epsilon(\mathbf{x};\mu) \mathbf{E}(\mathbf{x};\mu)) &=~ \rho(\mathbf{x};\mu), \quad \phantom{xxxi} (\mathbf{x};\mu) \in \Omega\times \mathcal{P},\\
\mathbf{E}(\mathbf{x};\mu)\times \vec{\mathbf{n}}(\mathbf{x}) &=~ \mathbf{0}, \quad \phantom{xxxxxxxx} (\mathbf{x};\mu) \in \partial\Omega\times \mathcal{P},
\end{cases}
\end{equation}
where $\vec{\mathbf{n}} := \vec{\mathbf{n}}(\mathbf{x})$ is the unit outward normal on $\partial\Omega$. In \cref{9-2-18ct1} the dielectric $\epsilon :=\epsilon (\mathbf{x};\mu)$, the magnetic permeability $\sigma :=\sigma(\mathbf{x};\mu)$ and the charge density $\rho :=\rho(\mathbf{x};\mu)$ are assumed to be known with 
\begin{align}\label{11-2-18ct5}
\underline{\rho} \le\rho(\mathbf{x};\mu) \le \overline{\rho}, \quad \underline{\epsilon} \le \epsilon (\mathbf{x};\mu) \le \overline{\epsilon} \quad\mbox{and}\quad \underline{\sigma} \le \sigma (\mathbf{x};\mu) \le \overline{\sigma}
\end{align}
a.e. in $\mathbf{x} \in \Omega$, all $\mu\in \mathcal{P}$ for some given constants $\underline{\rho}, \overline{\rho}, \underline{\epsilon}, \overline{\epsilon}, \underline{\sigma}$, and $\overline{\sigma}$ independent of both $\mathbf{x}$ and $\mu$, where $ \underline{\epsilon} >0 $ and $\underline{\sigma} >0$. Furthermore, we assume that Gauss's law is applied to the current source, i.e.
\begin{equation}\label{9-2-18ct2}
\nabla\cdot \left( \epsilon(\mathbf{x};\mu) \mathbf{u}(\mathbf{x})\right)  =0, \quad (\mathbf{x};\mu) \in \Omega\times\mathcal{P}.
\end{equation}
The function $\mathbf{u}$ denotes the control variable. For 
%For the system \cref{9-2-18ct1} if
\begin{equation*}
\mathbf{u} \in \mathbf{U}_{ad} := \left\{ \mathbf{u} \in \mathbf{L}^2(\Omega) := \left(L^2(\Omega)\right)^3 ~\big|~ \nabla\cdot (\epsilon \mathbf{u}) =0 \quad\mbox{and}\quad \underline{\mathbf{u}} \le \mathbf{u}\le \overline{\mathbf{u}}\right\}
\end{equation*}  
given, we solve \eqref{9-2-18ct1} for the electric field $\mathbf{E} := \mathbf{E}(\mathbf{x},\mathbf{u};\mu) := \mathbf{E}(\mathbf{u};\mu)\in \mathbf{H}_0(\mathbf{curl};\Omega)$ depending on $\mathbf{u}$ and the parameter $\mu$ as well (see \Cref{12-6-19ct2} for the definition of functional spaces). 
Here $\underline{\mathbf{u}}$ and  $\overline{\mathbf{u}} \in \mathbb{R}^3$ are given lower and upper bounds of the control. Therefore, for any given $\mu\in\mathcal{P}$ we define the {\it control-to-state} operator $\mathbf{E} : \mathbf{U}_{ad} \to \mathbf{H}_0(\mathbf{curl};\Omega)$ that maps each $\mathbf{u}$ to the unique weak solution $\mathbf{E}(\mb{u};\mu) $ of \cref{9-2-18ct1}. 

Let $D$ be a measurable subset of $\Omega$ and $\mathbf{E}_d(\mu) \in \mathbf{L}^2(D)$, $\mathbf{u}_d(\mu) \in \mathbf{L}^2(\Omega)$ respectively be the desired state and control,  both of which can be parameter dependent. 
In this paper we consider the parameterized control problem
$$
\min_{(\mathbf{u},\mathbf{E}) \in \mathbf{U}_{ad}\times \mathbf{H}_0(\mathbf{curl};\Omega)} J(\mathbf{u},\mathbf{E};\mu), \eqno(\mathbb{P}_e)
$$
where\footnote{The subscript $e$ in the problem $(\mathbb{P}_e)$ refers to ``exact".} the cost functional is defined as
$$
J(\mathbf{u},\mathbf{E};\mu) := \frac{1}{2}\|\sqrt{\epsilon(\mu)}(\mathbf{E}(\mu) - \mathbf{E}_d(\mu))\|^2_{\mathbf{L}^2(D)} + \frac{\alpha}{2} \|\sqrt{\epsilon(\mu)}(\mathbf{u} - \mathbf{u}_d(\mu))\|^2_{\mathbf{L}^2(\Omega)}
$$
and $\alpha>0$ is the regularization parameter. We assume that the desired state and control are uniformly $\mb{L}^2$-bounded with respect to the parameter, i.e. 
\begin{align}\label{8-5-19ct10}
\|\mb{E}_d(\mu)\|_{\mb{L}^2(D)} \le e_d \quad\mbox{and}\quad \|\mb{u}_d(\mu)\|_{\mb{L}^2(\Omega)} \le u_d
\end{align}
for all $\mu\in\mathcal{P}$ with $e_d$ and $u_d$ some positive constants. Furthermore, $\mathbf{E}_d$ fulfills the Gauss's law in $D$, i.e., 
\begin{align*}
\nabla\cdot (\epsilon(\mu) \mathbf{E}_d(\mu)) = \rho(\mu) \quad \mbox{in } D . 
\end{align*}
%is satisfied on $D$. 

Let $(\mathcal{E}_h,V_h)$ be the finite element space associated with the system \cref{9-2-18ct1} and $\mathbf{E}_h(\mu)$ be the finite element approximation of $\mathbf{E}(\mu)$ (cf. \Cref{fem}). Adopting the variational discretization concept introduced in \cite{hin05} (where control is not directly discretized), we approximate the ``exact" problem $(\mathbb{P}_e)$ by the discrete one
$$\min_{(\mathbf{u},\mathbf{E}_h) \in \mathbf{U}_{ad}\times \mathcal{E}_h} J(\mathbf{u},\mathbf{E}_h;\mu), \eqno(\mathbb{P}_h)$$
%where 
%$$J(\mathbf{u},\mathbf{E}_h;\mu) := \frac{1}{2}\|\sqrt{\epsilon(\mu)}(\mathbf{E}_h(\mu) - \mathbf{E}_d(\mu))\|^2_{\mathbf{L}^2(D)} + \frac{\alpha}{2} \|\sqrt{\epsilon(\mu)}(\mathbf{u} - \mathbf{u}_d(\mu))\|^2_{\mathbf{L}^2(\Omega)}$$
subject to
\begin{equation}\label{21-8-19ct1}
\begin{cases}
\left(\sigma^{-1}(\mu)\nabla\times \mathbf{E}_h(\mu), \nabla\times \mathbf{\Phi}_h\right)_{\mathbf{L}^2(\Omega)} &=~ \left(\epsilon(\mu) \mathbf{u}, \mathbf{\Phi}_h\right)_{\mathbf{L}^2(\Omega)}\\
\left(\epsilon(\mu) \mathbf{E}_h(\mu), \nabla\phi_h\right)_{\mathbf{L}^2(\Omega)} &=~ -(\rho(\mu),\phi_h)_{L^2(\Omega)}
\end{cases}
\end{equation}
for all $(\mathbf{\Phi}_h,\phi_h) \in\mathcal{E}_h \times V_h$.

As mentioned in the introduction, the constrained saddle point system \eqref{21-8-19ct1} is expensive to solve. Our goal is to create a reduced basis method for $(\mathbb{P}_h)$, prove its convergence, and derive a posteriori error estimates.

\section{Preliminaries}\label{12-6-19ct2}

We start this section by presenting the definition of functional spaces which are utilized in the paper, for more details one can consult \cite{amrouche,monk}. Well-posedness and finite element discretization, including a priori error estimates, of \eqref{9-2-18ct1} are given in Subsections~\ref{9-2-18ct1} and \ref{fem}, respectively.

\subsection{Functional spaces}

In this paper bold typeface is used to indicate a point in $\mathbb{R}^3$, a (three-dimensional) vector-valued function or a Hilbert space of vector-valued functions.
% As usual, we use the differential operator $\nabla := \left(\dfrac{\partial}{\partial x^1}, \dfrac{\partial}{\partial x^2}, \dfrac{\partial}{\partial x^3}\right)$\HA{, so that} $\nabla\phi(\mathbf{x}) = \left(\dfrac{\partial\phi(\mathbf{x})}{\partial x^1}, \dfrac{\partial\phi(\mathbf{x})}{\partial x^2}, \dfrac{\partial\phi(\mathbf{x})}{\partial x^3}\right)$ for the scalar function $\phi = \phi(\mathbf{x})$. Furthermore, for the vector-valued function $\mathbf{\Phi} = \mathbf{\Phi}(\mathbf{x}) = \left(\Phi^1(\mathbf{x}), \Phi^2(\mathbf{x}), \Phi^3(\mathbf{x})\right)$ we have
%\begin{align*}
%&\nabla \cdot \mathbf{\Phi} = \dfrac{\partial\Phi^1(\mathbf{x})}{\partial x^1} + \dfrac{\partial \Phi^2(\mathbf{x})}{\partial x^2} + \dfrac{\partial\Phi^3(\mathbf{x})}{\partial x^3} := \mbox{div} \mathbf{\Phi} \quad\mbox{and}\\
%&\nabla \times \mathbf{\Phi} = \left( \dfrac{\partial\Phi^3(\mathbf{x})}{\partial x^2} - \dfrac{\partial\Phi^2(\mathbf{x})}{\partial x^3}, \dfrac{\partial\Phi^1(\mathbf{x})}{\partial x^3} - \dfrac{\partial\Phi^3(\mathbf{x})}{\partial x^1}, \dfrac{\partial\Phi^2(\mathbf{x})}{\partial x^1} - \dfrac{\partial\Phi^1(\mathbf{x})}{\partial x^2}\right) := \mathbf{curl} \mathbf{\Phi}.
%\end{align*}
%We also mention that $\nabla\times (\nabla\phi) = \mathbf{0}$ and $\nabla\cdot (\nabla\times\mathbf{\Phi}) =0$. 
The Hilbert spaces
\begin{align*}
\mathbf{H}(\mbox{div}; \Omega) 
&:= \left\{ \mathbf{\Phi} \in \mathbf{L}^2(\Omega) ~\big|~ \nabla\cdot\mathbf{\Phi} \in L^2(\Omega)\right\} \quad\mbox{and}\\
\mathbf{H}(\mathbf{curl};\Omega) 
&:= \left\{ \mathbf{\Phi} \in \mathbf{L}^2(\Omega) ~\big|~ \nabla\times\mathbf{\Phi} \in \mathbf{L}^2(\Omega)\right\}
\end{align*}
are respectively equipped the inner product
\begin{align*}
(\mathbf{\Phi}, \mathbf{\Psi})_{\mathbf{H}(\mbox{div}; \Omega)} 
&:= (\mathbf{\Phi}, \mathbf{\Psi})_{\mathbf{L}^2(\Omega)} + (\nabla\cdot \mathbf{\Phi}, \nabla\cdot \mathbf{\Psi})_{L^2(\Omega)} \quad\mbox{and}\\
(\mathbf{\Phi}, \mathbf{\Psi})_{\mathbf{H}(\mathbf{curl}; \Omega)} 
&:= (\mathbf{\Phi}, \mathbf{\Psi})_{\mathbf{L}^2(\Omega)} + (\nabla\times \mathbf{\Phi}, \nabla\times \mathbf{\Psi})_{\mathbf{L}^2(\Omega)}.
\end{align*}
%\HA{Notice that}
%$$\mathbf{H}(\mbox{div}; \Omega) = {\overline{\mathbf{C}^\infty(\overline{\Omega})}}^{\mathbf{H}(\mbox{div}; \Omega)} \quad\mbox{and}\quad \mathbf{H}(\mathbf{curl};\Omega) = {\overline{\mathbf{C}^\infty(\overline{\Omega})}}^{\mathbf{H}(\mathbf{curl};\Omega)},$$
%where the closures are taken with respect to the norm of the space $\mathbf{H}(\mbox{div}; \Omega)$ and $\mathbf{H}(\mathbf{curl};\Omega)$, respectively.

The {\it normal trace operator} $\gamma_n(\mathbf{\Phi}) := \vec{\mathbf{n}} \cdot \mathbf{\Phi}_{|\partial\Omega} $ for all $\mathbf{\Phi}\in \mathbf{C}^\infty(\overline{\Omega})$ can be extended to a surjective, continuous linear map from $\mathbf{H}(\mbox{div}; \Omega) \to H^{-1/2}(\partial\Omega) := \left(H^{1/2}(\partial\Omega)\right)^*$ such that Green's formula (cf.\ \cite[\S 3]{monk})
\begin{align}\label{11-2-18ct1}
(\nabla\cdot \mathbf{\Phi}, \phi)_{L^2(\Omega)} = -(\mathbf{\Phi}, \nabla\phi)_{\mathbf{L}^2(\Omega)} + \left\langle \gamma_n(\mathbf{\Phi}), \phi \right\rangle_{\left(H^{-1/2}(\partial\Omega), H^{1/2}(\partial\Omega)\right)}
\end{align}
holds true for all $\mathbf{\Phi} \in \mathbf{H}(\mbox{div}; \Omega)$ and $\phi\in H^1(\Omega)$. The {\it tangential trace operator} $\gamma_t(\mathbf{\Phi}) :=  \vec{\mathbf{n}} \times \mathbf{\Phi}_{|\partial\Omega}$ for all $\mathbf{\Phi}\in \mathbf{C}^\infty(\overline{\Omega})$ can be also extended to a continuous linear map from $\mathbf{H}(\mathbf{curl}; \Omega) \to \mathbf{H}^{-1/2}(\partial\Omega)$. Further, Green's formula \cite[Theorem 3.29]{monk} 
\begin{align}\label{11-2-18ct3}
(\nabla\times \mathbf{\Phi}, \mathbf{\Psi})_{\mathbf{L}^2(\Omega)} = (\mathbf{\Phi}, \nabla\times\mathbf{\Psi})_{\mathbf{L}^2(\Omega)} + \left\langle \gamma_t(\mathbf{\Phi}), \mathbf{\Psi} \right\rangle_{\left(\mathbf{H}^{-1/2}(\partial\Omega), \mathbf{H}^{1/2}(\partial\Omega)\right)}
\end{align}
holds true for all $\mathbf{\Phi} \in \mathbf{H}(\mathbf{curl}; \Omega)$ and $\mathbf{\Psi}\in \mathbf{H}^1(\Omega)$.

We conclude this subsection by the following definition
$$\mathbf{H}_0(\mbox{div}; \Omega) :=  {\overline{\mathbf{C}^\infty_c(\Omega)}}^{\mathbf{H}(\mbox{div}; \Omega)} \quad\mbox{and}\quad \mathbf{H}_0(\mathbf{curl};\Omega) := {\overline{\mathbf{C}^\infty_c(\Omega)}}^{\mathbf{H}(\mathbf{curl};\Omega)},$$
where the closures are respectively taken with respect to the norm of the space $\mathbf{H}(\mbox{div}; \Omega)$ and $\mathbf{H}(\mathbf{curl};\Omega)$ and  $C^\infty_c(\Omega)$ is the space of all infinitely continuously differential functions with compact support in $\Omega$. Notice that
\begin{align*}
\mathbf{H}_0(\mbox{div}; \Omega) &:=  \left\{ \mathbf{\Phi} \in \mathbf{H}(\mbox{div};\Omega) ~\big|~ \gamma_n(\mathbf{\Phi}) =0 \right\}, \\
\mathbf{H}_0(\mathbf{curl};\Omega) &:= \left\{ \mathbf{\Phi} \in \mathbf{H}(\mathbf{curl};\Omega) ~\big|~  \gamma_t(\mathbf{\Phi}) = \mathbf{0}\right\}.
\end{align*}

\subsection{Variational formulation of the system \cref{9-2-18ct1}}

For any given $\mu\in\mathcal{P}$ and $\mb{u}\in\mb{L}^2(\Omega)$ an element $\mathbf{E} := \mb{E}(\mu) := \mb{E}(\mb{u};\mu)\in \mathbf{H}_0(\mathbf{curl};\Omega)$ is said to be a weak solution of \cref{9-2-18ct1} if 
\begin{equation}\label{11-2-18ct4}
\begin{cases}
\left(\sigma^{-1}(\mu)\nabla\times \mathbf{E}(\mu), \nabla\times \mathbf{\Phi}\right)_{\mathbf{L}^2(\Omega)} &=~ \left(\epsilon(\mu) \mathbf{u}, \mathbf{\Phi}\right)_{\mathbf{L}^2(\Omega)}, \quad\forall \mathbf{\Phi}\in\mathbf{H}_0(\mathbf{curl};\Omega)\\
\left(\epsilon(\mu) \mathbf{E}(\mu), \nabla\phi\right)_{\mathbf{L}^2(\Omega)} &=~ -(\rho(\mu),\phi)_{L^2(\Omega)},~ \quad\forall \phi\in H^1_0(\Omega).
\end{cases}
\end{equation}
The first equation in \eqref{11-2-18ct4} is obtained by multiplying the first equation of \cref{9-2-18ct1} with $\mathbf{\Phi}\in\mathbf{H}_0(\mathbf{curl};\Omega)$, and then using the identity \cref{11-2-18ct3}. The second equation of \cref{11-2-18ct4} is obtained by using the second equation in \cref{9-2-18ct1} and the Green's formula \cref{11-2-18ct1}. 

We define
\begin{align*}
\mathbf{V} 
&:= \left\{ \mathbf{\tau}\in \mathbf{H}_0(\mathbf{curl};\Omega) ~|~ \nabla\cdot(\epsilon\mathbf{\mathbf{\tau}}) =0 \right\} .
\end{align*}
Then by the compactness of the embedding $\mathbf{V} \hookrightarrow \mathbf{L}^2(\Omega)$ (for $\epsilon$ piecewise smooth) and 
$$
\left\{\mathbf{\tau} \in \mathbf{V} ~|~ \nabla \times \mathbf{\tau} = \mathbf{0} \right\}
= \left\{ \mathbf{\tau} \in \mathbf{L}^2(\Omega) ~\big|~ \nabla\times\mathbf{\tau} = \mathbf{0},~ \nabla \cdot (\epsilon\mathbf{\tau}) = 0,~ \vec{\mathbf{n}} \times \mathbf{\tau}_{|\partial\Omega} = \mathbf{0}\right\} = \{\mathbf{0}\}
$$
(see, \cite{cos90,yousept-zou}), an application of Peetre's lemma (see, \cite[Lemma~2]{kangro}) yields that there exists a positive $C^\Omega$ such that
\begin{align}\label{12-6-18ct1}
\|\mathbf{\tau}\|_{\mathbf{L}^2(\Omega)} \le C^\Omega\| \nabla\times \mathbf{\tau}\|_{\mathbf{L}^2(\Omega)}
\end{align}
for all $\mathbf{\tau}\in \mathbf{V}$. Thus, with the aid of the condition \cref{11-2-18ct5}, we have the coercivity condition
\begin{align}\label{23-6-19ct2}
\|\mathbf{v}\|^2_{\mathbf{H}(\mathbf{curl};\Omega)} \le C^\Omega_{\overline{\sigma}} \left(\sigma^{-1}\nabla\times \mathbf{v}, \nabla\times \mathbf{v}\right)_{\mathbf{L}^2(\Omega)}
\end{align}
for all $\mathbf{v}\in \mathbf{V}$, where the constant $C^\Omega_{\overline{\sigma}} >0$ is independent of $\mathbf{v}$. Due to the standard theory of the mixed variational problems  (see, e.g., \cite{brezzi,gatica}), we conclude that the system \cref{11-2-18ct4}  attains a unique solution  $\mathbf{E} = \mathbf{E}(\mb{u};\mu)\in \mathbf{H}_0(\mathbf{curl};\Omega)$ which satisfies
\begin{align}\label{3-5-19ct7}
\|\mathbf{E}\|_{\mathbf{H}(\mathbf{curl};\Omega)} \le C^\Omega_{\overline{\epsilon},\overline{\sigma}}\left(\| \mathbf{u}\|_{\mathbf{L}^2(\Omega)} + \|\rho\|_{L^2(\Omega)}\right) \le C^\Omega_{\overline{\epsilon},\overline{\sigma}}|\Omega|^{1/2}(\overline{\rho} +\|\overline{\mb{u}}\|_{\mathbb{R}^3}) := C_\mathbf{E}.
\end{align}
We therefore can define for any fixed $\mu\in\mathcal{P}$ the map
\begin{align*}
S_\mu : \mathbf{U}_{ad} \rightarrow \mathbf{H}_0(\mathbf{curl};\Omega) \quad\mbox{with} \quad \mathbf{u} \mapsto S_\mu(\mb{u}) := \mathbf{E}(\mb{u})
\end{align*}
and for any fixed $\mb{u} \in \mathbf{U}_{ad}$ the one
\begin{align*}
S_{\mb{u}} :  \mathcal{P} \rightarrow \mathbf{H}_0(\mathbf{curl};\Omega) \quad\mbox{with} \quad \mu \mapsto S_{\mb{u}}(\mu) := \mathbf{E}(\mu).
\end{align*}

\subsection{Finite element discretization}\label{fem}

Hereafter, we assume $\left(\mathcal{T}_h\right)_{h>0}$ is a quasi-uniform family of regular triangulations  of the domain $\Omega$ with the mesh size $h$
(cf.\ \cite{Brenner_Scott}). For discretization of the state variable solving the system \cref{11-2-18ct4} let us denote the N\'ed\'elec finite element spaces (cf.\ \cite{nedelec}) 
\begin{equation*}
\begin{aligned}
&\mathcal{E}_h := \left\{\mathbf{E}_h \in \mathbf{H}_0(\mathbf{curl};\Omega) ~\big|~ {\mathbf{E}_h}_{|T} = \mathbf{a}_T + \mathbf{b}_T \times \mathbf{x}, \enskip\forall T\in \mathcal{T}_h \enskip\mbox{with}\enskip \mathbf{a}_T, \mathbf{b}_T \in \mathbb{R}^3\right\},\\
& V_h := \left\{\phi_h \in H^1_0(\Omega) ~\big|~ {\phi_h}_{|T} = a_T + \mathbf{b}_T \cdot \mathbf{x}, \enskip\forall T\in \mathcal{T}_h \enskip\mbox{with}\enskip a_T\in \mathbb{R}, \mathbf{b}_T \in \mathbb{R}^3\right\},
\end{aligned}
\end{equation*}
where $\nabla V_h \subset \mathcal{E}_h$.

The discrete variational formulation corresponding to the system \cref{11-2-18ct4} then reads: find $\mathbf{E}_h \in \mathcal{E}_h$ such that \cref{21-8-19ct1} is satisfied for all $(\mathbf{\Phi}_h,\phi_h) \in\mathcal{E}_h \times V_h$.
Similarly to \cref{23-6-19ct2}, since the discrete Poincar\'e-Friedrichs-type inequality (cf.\ \cite[Theorem 4.7]{hip02}, \cite[\S 7]{monk})
\begin{align}\label{13-2-18ct7}
\|\mathbf{v}_h\|_{\mathbf{L}^2(\Omega)} \le C^\Omega\|\nabla\times \mathbf{v}_h\|_{\mathbf{L}^2(\Omega)}
\end{align}
is satisfied for all discrete $\epsilon$-divergence-free functions, i.e.
\begin{align}\label{eq:Dheps}
\mathbf{v}_h \in \mathcal{D}_h^{(\epsilon)} := \left\{ \mathbf{E}_h \in \mathcal{E}_h ~\big|~ \left(\epsilon \mathbf{E}_h, \nabla\phi_h\right)_{\mathbf{L}^2(\Omega)} = 0 \quad\mbox{for all}\quad \phi_h\in V_h\right\},
\end{align}
we have that 
\begin{align}\label{30-4-19ct1}
\|\mathbf{v}_h\|^2_{\mathbf{H}(\mathbf{curl};\Omega)} \le C^\Omega_{\overline{\sigma}} (\sigma^{-1}\nabla\times \mathbf{v}_h, \nabla\times \mathbf{v}_h)_{\mathbf{L}^2(\Omega)}
\end{align}
for all $\mathbf{v}_h\in \mathcal{D}_h^{(\epsilon)}$. Therefore, we conclude that the system \cref{21-8-19ct1} has a unique solution  $\mathbf{E}_h\in \mathcal{E}_h$ satisfying the estimate
$\|\mathbf{E}_h\|_{\mathbf{H}(\mathbf{curl};\Omega)} \le C_\mathbf{E}.
$

%Furthermore, it seems that the numerical computation for \cref{16-2-18ct1} may be easier than for the system \cref{21-8-19ct1}. We in the current work thus consider $\tilde{\mathbf{E}}_h$ from \cref{16-2-18ct1} as an approximation of the solution $\mathbf{E}\in \mathbf{H}_0(\mathbf{curl};\Omega)$ to the system \cref{11-2-18ct4}.

For all $s\ge0$ we denote by
$$\mathbf{H}^s(\mathbf{curl};\Omega) := \left\{ \mathbf{\Phi} \in \mathbf{H}^s(\Omega) ~\big|~ \nabla\times \mathbf{\Phi} \in \mathbf{H}^s(\Omega)\right\}.$$
Equipped with the norm
$$\|\mathbf{\Phi}\|_{\mathbf{H}^s(\mathbf{curl};\Omega)} := \left( \|\mathbf{\Phi}\|^2_{\mathbf{H}^s(\Omega)} + \|\nabla\times \mathbf{\Phi}\|^2_{\mathbf{H}^s(\Omega)} \right)^{1/2},$$
it is a Banach space. Before going further, we state the following result. 

\begin{theorem}\label{23-8-19ct1}
For any given $\mu\in\mathcal{P}$ let $\mathbf{E}(\mu)$ and $\mathbf{E}_h(\mu)$ be the unique solution to \cref{11-2-18ct4} and \cref{21-8-19ct1}, respectively. Then:

(i) There holds the limit
\begin{align*}
\lim_{h\to 0} \|\nabla \times (\mathbf{E}(\mu) - \mathbf{E}_h(\mu))\|_{\mathbf{L}^2(\Omega)} =0.
\end{align*}

(ii) In addition $\epsilon(\mu), \sigma^{-1}(\mu) \in W^{1,\infty}(\Omega)$ we get the regularity $\mathbf{E}(\mu) \in \mathbf{H}^s(\mathbf{curl};\Omega)$ for some $s\in (1/2, 1]$. Furthermore, there exist constants $\nu, \nu'\in (1/2, 1]$ such that the estimates
\begin{align*}
\|\mathbf{E}(\mu) - \mathbf{E}_h(\mu)\|_{\mathbf{H}(\mathbf{curl};\Omega)}
&\le Ch^s \|\mathbf{E}(\mu)\|_{\mathbf{H}^s(\mathbf{curl};\Omega)}\\
\|\nabla \cdot \epsilon(\mathbf{E}(\mu) - \mathbf{E}_h(\mu))\|_{H^{-\nu}(\Omega)}
&\le Ch^{\nu+\nu'-1} \left(\|\nabla\times \mb{E}(\mu)\|_{\mathbf{L}^2(\Omega)} + \|\rho(\mu)\|_{H^{\nu'-1}(\Omega)}\right)
\end{align*}
hold true.
\end{theorem}

\begin{proof}
The regularity $\mathbf{E} \in \mathbf{H}^s(\mathbf{curl};\Omega)$ follows from \cite[Lemma 3.6]{ciarlet-wu-zou}. Further, the assertion is based on standard arguments, it is therefore omitted here.
\end{proof}

\section{Primal reduced basis approach}\label{RBM}

By standard arguments (see, e.g., \cite{hin08,tro10}), one can verify that the problem $(\mathbb{P}_e)$ attains a unique solution for each the parameter $\mu \in \mathcal{P}$. Furthermore, we can derive the following, for instance using Lagrangian approach, first order optimality system satisfied by the optimal control,  state and  adjoint. 

\begin{theorem}\label{24-4-19ct2}
Given $\mu\in\mathcal{P}$, the pair $(\mathbf{u}^*_e(\mu),\mathbf{E}^*_e(\mu)) \in \mathbf{U}_{ad}\times \mathbf{H}_0(\mathbf{curl};\Omega)$ is the unique solution\footnote{The superscript $^*$ refers to ``optimality".} of the problem $(\mathbb{P}_e)$ if and only if there exists an adjoint state  $\mathbf{F}^*_e(\mu) \in \mathbf{H}_0(\mathbf{curl};\Omega)$ such that the triple $(\mathbf{u}^*_e(\mu),\mathbf{E}^*_e(\mu),\mathbf{F}^*_e(\mu))$ satisfies the system 
\begin{subequations}
\begin{align}
&\left(\sigma^{-1}(\mu)\nabla\times \mathbf{E}^*_e(\mu), \nabla\times \mathbf{\Phi}\right)_{\mathbf{L}^2(\Omega)} = \left(\epsilon(\mu) \mathbf{u}^*_e(\mu), \mathbf{\Phi}\right)_{\mathbf{L}^2(\Omega)}, \label{24-4-19ct3}\\
&\left(\epsilon(\mu) \mathbf{E}^*_e(\mu), \nabla\phi\right)_{\mathbf{L}^2(\Omega)} = -(\rho(\mu),\phi)_{L^2(\Omega)}, \label{24-4-19ct4}\\
&\left(\sigma^{-1}(\mu)\nabla\times \mathbf{F}^*_e(\mu), \nabla\times \mathbf{\Phi}\right)_{\mathbf{L}^2(\Omega)} = \left(\epsilon(\mu) (\mathbf{E}^*_e(\mu) - \mathbf{E}_d(\mu)), \mathbf{\Phi}\right)_{\mathbf{L}^2(D)}, \label{24-4-19ct5}\\
&\left(\epsilon(\mu) \mathbf{F}^*_e(\mu), \nabla\phi\right)_{\mathbf{L}^2(\Omega)} = 0, \label{24-4-19ct6}\\
&\left(\epsilon(\mu)\big(\mathbf{u} - \mathbf{u}^*_e(\mu)\big), \mathbf{u}_d(\mu)  - \frac{1}{\alpha}\mathbf{F}^*_e(\mu) - \mathbf{u}^*_e(\mu)\right)_{\mathbf{L}^2(\Omega)} \le 0 \label{24-4-19ct7}
\end{align}
\end{subequations}
for all $\left(\mb{\Phi},\phi,\mb{u} \right) \in \mathbf{H}_0(\mathbf{curl};\Omega) \times H^1_0(\Omega) \times \mb{U}_{ad}$. 
\end{theorem}
%Since $\nabla \cdot (\epsilon \mb{u}^*_e) = \nabla \cdot \left( \mb{E}^*_e - \mb{E}_d\right) = 0$, the above system is well defined.
Notice that the following inequality holds true for all $\mu\in\mathcal{P}$
\begin{align}\label{12-4-20ct1}
\|\mathbf{F}^*_e(\mu)\|_{\mathbf{L}^2(\Omega)} \le
C^\Omega_{\overline{\epsilon},\overline{\sigma}} \left( e_d + C_\mathbf{E}\right) := C_\mathbf{F}
\end{align}
Based on the finite element approach in \Cref{fem}, next we approximate the ``exact" problem $(\mathbb{P}_e)$ by the discrete one $(\mathbb{P}_h)$. Then, the associated first order optimality system for the problem $(\mathbb{P}_h)$ reads:

\begin{theorem}\label{24-4-19ct2*}
Given $\mu\in\mathcal{P}$, the pair $(\mathbf{u}^*_h(\mu),\mathbf{E}^*_h(\mu)) \in \mathbf{U}_{ad}\times \mathcal{E}_h$ is the unique solution of the problem $(\mathbb{P}_h)$ if and only if there exists an adjoint state  $\mathbf{F}^*_h(\mu) \in \mathcal{E}_h$ such that the triple $(\mathbf{u}^*_h(\mu),\mathbf{E}^*_h(\mu),\mathbf{F}^*_h(\mu))$ satisfies the system
\begin{subequations}
\begin{align}
&\left(\sigma^{-1}(\mu)\nabla\times \mathbf{E}^*_h(\mu), \nabla\times \mathbf{\Phi}_h\right)_{\mathbf{L}^2(\Omega)} = \left(\epsilon(\mu) \mathbf{u}^*_h(\mu), \mathbf{\Phi}_h\right)_{\mathbf{L}^2(\Omega)}, \label{24-4-19ct3*}\\
&\left(\epsilon(\mu) \mathbf{E}^*_h(\mu), \nabla\phi_h\right)_{\mathbf{L}^2(\Omega)} = -(\rho(\mu),\phi_h)_{L^2(\Omega)}, \label{24-4-19ct4*}\\
&\left(\sigma^{-1}(\mu)\nabla\times \mathbf{F}^*_h(\mu), \nabla\times \mathbf{\Phi}_h\right)_{\mathbf{L}^2(\Omega)} = \left(\epsilon(\mu) (\mathbf{E}^*_h(\mu) - \mathbf{E}_d(\mu)), \mathbf{\Phi}_h\right)_{\mathbf{L}^2(D)}, \label{24-4-19ct5*}\\
&\left(\epsilon(\mu) \mathbf{F}^*_h(\mu), \nabla\phi_h\right)_{\mathbf{L}^2(\Omega)} = 0, \label{24-4-19ct6*}\\
&\left(\epsilon(\mu)\big(\mathbf{u} - \mathbf{u}^*_h(\mu)\big), \mathbf{u}_d(\mu)  - \frac{1}{\alpha}\mathbf{F}^*_h(\mu) - \mathbf{u}^*_h(\mu)\right)_{\mathbf{L}^2(\Omega)} \le 0 \label{24-4-19ct7*}
\end{align}
\end{subequations}
for all $\left( \mb{\Phi}_h, \phi_h, \mb{u}\right) \in \mathcal{E}_h \times V_h \times \mb{U}_{ad}$.
\end{theorem}

The above optimality system \cref{24-4-19ct3*}--\cref{24-4-19ct7*} constitutes several sets of variational equations and  inequalities 
which may be computationally expensive. Thus the surrogate model approach will be considered next, where the original high 
dimensional problem  is replaced by a reduced order approximation.
%
%
% to solve when the dimension of the finite element spaces is  very high, and in particular the variety of parameters is taken into account. From a reduced basis point of view, one needs a surrogate model for the system which attains two purposes. It first approximates well enough for the original system and second has a significantly lower dimension. To do so, a crucial requirement is a careful choice of these parameters (cf.\ \cite{Ha17,KaGr14,Ne13,To11}). 

Assume that we are given the reduced basis spaces
$$(\mathcal{E}_N,V_N) \subset (\mathcal{E}_h,V_h).$$
%that are orthonormal  due to, for example, a standard Gram-Schmidt orthogonalisation process. 
Furthermore, to guarantee the existence of a solution to the constraint system, we assume that the coercivity condition \cref{30-4-19ct1} is fulfilled on $(\mathcal{E}_N,V_N)$.
We can then consider the reduced basis problem 
%$(\mathbb{P}_N)$ corresponding to the full approximation $(\mathbb{P}_h)$, i.e.
$$\min_{(\mathbf{u},\mathbf{E}_N) \in \mathbf{U}_{ad}\times \mathcal{E}_N} J(\mathbf{u},\mathbf{E}_N;\mu), \eqno(\mathbb{P}_N)$$
%where
%$$J(\mathbf{u},\mathbf{E}_N;\mu) := \frac{1}{2}\|\sqrt{\epsilon(\mu)}(\mathbf{E}_N(\mu) - \mathbf{E}_d(\mu))\|^2_{\mathbf{L}^2(D)} + \frac{\alpha}{2} \|\sqrt{\epsilon(\mu)}(\mathbf{u} - \mathbf{u}_d(\mu))\|^2_{\mathbf{L}^2(\Omega)}$$
subject to
\begin{equation}\label{21-8-19ct2}
\begin{cases}
\left(\sigma^{-1}(\mu)\nabla\times \mathbf{E}_N(\mu), \nabla\times \mathbf{\Phi}_N\right)_{\mathbf{L}^2(\Omega)} &=~ \left(\epsilon(\mu) \mathbf{u}, \mathbf{\Phi}_N\right)_{\mathbf{L}^2(\Omega)}\\
\left(\epsilon(\mu) \mathbf{E}_N(\mu), \nabla\phi_N\right)_{\mathbf{L}^2(\Omega)} &=~ -(\rho(\mu),\phi_N)_{L^2(\Omega)}
\end{cases}
\end{equation}
for all $(\mathbf{\Phi}_N,\phi_N) \in\mathcal{E}_N \times V_N$. 
The associated first order optimality system reads:

\begin{theorem}\label{24-4-19ct2**}
Given $\mu\in\mathcal{P}$, the pair $(\mathbf{u}^*_N(\mu),\mathbf{E}^*_N(\mu)) \in \mathbf{U}_{ad}\times \mathcal{E}_N$ is the unique solution of the problem $(\mathbb{P}_N)$ if and only if there exists an adjoint state  $\mathbf{F}^*_N(\mu) \in \mathcal{E}_N$ such that the triple $(\mathbf{u}^*_N(\mu),\mathbf{E}^*_N(\mu),\mathbf{F}^*_N(\mu))$ satisfies the system
\begin{subequations}
\begin{align}
&\left(\sigma^{-1}(\mu)\nabla\times \mathbf{E}^*_N(\mu), \nabla\times \mathbf{\Phi}_N\right)_{\mathbf{L}^2(\Omega)} = \left(\epsilon(\mu) \mathbf{u}^*_N(\mu), \mathbf{\Phi}_N\right)_{\mathbf{L}^2(\Omega)}, \label{24-4-19ct3**}\\
&\left(\epsilon(\mu) \mathbf{E}^*_N(\mu), \nabla\phi_N\right)_{\mathbf{L}^2(\Omega)} = -(\rho(\mu),\phi_N)_{L^2(\Omega)}, \label{24-4-19ct4**}\\
&\left(\sigma^{-1}(\mu)\nabla\times \mathbf{F}^*_N(\mu), \nabla\times \mathbf{\Phi}_N\right)_{\mathbf{L}^2(\Omega)} = \left(\epsilon(\mu) (\mathbf{E}^*_N(\mu) - \mathbf{E}_d(\mu)), \mathbf{\Phi}_N\right)_{\mathbf{L}^2(D)}, \label{24-4-19ct5**}\\
&\left(\epsilon(\mu) \mathbf{F}^*_N(\mu), \nabla\phi_N\right)_{\mathbf{L}^2(\Omega)} = 0, \label{24-4-19ct6**}\\
&\left(\epsilon(\mu)\big(\mathbf{u} - \mathbf{u}^*_N(\mu)\big), \mathbf{u}_d(\mu)  - \frac{1}{\alpha}\mathbf{F}^*_N(\mu) - \mathbf{u}^*_N(\mu)\right)_{\mathbf{L}^2(\Omega)} \le 0 \label{24-4-19ct7**}
\end{align}
\end{subequations}
for all $\left( \mb{\Phi}_N, \phi_N, \mb{u}\right) \in \mathcal{E}_N \times V_N \times \mb{U}_{ad}$.
\end{theorem}

We conclude this section by performing the greedy  sampling procedure \cite{HaOh08,Hesthaven16,AQ16} applied to the problem under consideration. Note that, by the discrete Helmholtz decomposition (see, \cite[\S 7.2.1]{monk}), for all $z^\epsilon_h \in \mathcal{E}_h$, there exists a unique pair $(z^1_h,H(z^\epsilon_h)) \in \mathcal{D}^{(1)}_h \times V_h$ such that $z^\epsilon_h = z^1_h + \nabla H(z^\epsilon_h)$. 

\begin{algorithm}
\caption{Greedy procedure}
\label{31-8-18ct1}
\begin{algorithmic}
\STATE{Choose $\mathcal{S}_{train} \subset \mathcal{P}$, an arbitrary $\mu^1 \in \mathcal{S}_{train}$, $\epsilon_{\mbox{tol}} >0$ and $N_{max} \in\mathbb{N}$}

\STATE{Set $N :=1$, and
\begin{align*} 
&\mathcal{E}_N := \mbox{span} \{\mathbf{E}^*_h(\mu^N),\mathbf{F}^*_h(\mu^N), \nabla H(\mathbf{E}^*_h(\mu^N)), \nabla H(\mathbf{F}^*_h(\mu^N))\}\\ 
&V_N := \mbox{span} \{H(\mathbf{E}^*_h(\mu^N)), H(\mathbf{F}^*_h(\mu^N))\}
\end{align*}
}
\WHILE{$\max_{\mu\in \mathcal{S}_{\mbox{train}}} \Delta_N(\mathcal{E}_N,V_N;\mu) > \epsilon_{\mbox{tol}}$ and $N\le N_{\mbox{max}}$}

\STATE{$\mu^{N+1} := \arg \max_{\mu\in \mathcal{S}_{\mbox{train}}} \Delta_N(\mathcal{E}_N,V_N;\mu)$}

\STATE{$\mathcal{E}_{N+1} := \mbox{span}\{\{\mathbf{E}^*_h(\mu^{N+1}),\mathbf{F}^*_h(\mu^{N+1}), \nabla H(\mathbf{E}^*_h(\mu^{N+1})), \nabla H(\mathbf{F}^*_h(\mu^{N+1})) \} \cup \mathcal{E}_N\}$}

\STATE{$V_{N+1} := \mbox{span}\{\{H(\mathbf{E}^*_h(\mu^{N+1})),H(\mathbf{F}^*_h(\mu^{N+1}))\} \cup V_N\}$}

\STATE{$N:= N+1$}

\ENDWHILE
\end{algorithmic}
\end{algorithm}

%\lipsum[41]
%\lipsum[50]

In \cref{31-8-18ct1}, the sampling parameter set $\mathcal{S}_{train} \subset \mathcal{P}$ is finite, but rich enough to so that it is a good approximation of the full parameter set $\mathcal{P}$. The initial parameter $\mu^1$ is chosen arbitrarily in $\mathcal{S}_{train}$, $\epsilon_{\mbox{tol}}$ is a desired error tolerance and $N_{max}$ is the maximum number of iterations. The pair $\{\mathbf{E}^*_h(\mu^N),\mathbf{F}^*_h(\mu^N)\}$ is the optimal state and adjoint state defined by the optimality system \cref{24-4-19ct3*}--\cref{24-4-19ct7*} at the parameter $\mu =\mu^N$. The quantity $\Delta_N(\mathcal{E}_N,V_N;\mu)$ is an error estimator between solutions of the problem $(\mathbb{P}_h)$ and the reduced one $(\mathbb{P}_N)$ at the given parameter $\mu$, that will be described the detail in \Cref{Post-estimate}.

\section{Convergence of the reduced basis method} \label{convergence}

Our aim in this section is to investigate the uniform convergence 
\begin{align*}
\lim_{N\to\infty}\sup_{\mu\in\mathcal{P}} \|\mathbf{u}^*_h(\mu) - \mathbf{u}^*_N(\mu)\|_{\mathbf{L}^2(\Omega)} =0
\end{align*}
of reduced basis optimal solutions to the original one. To do so, we assume that the snapshot parameter sample $\mathcal{P}_N := \{\mu^1, ..., \mu^N\}$, where $\mu_i$ are chosen via \Cref{31-8-18ct1}, is dense in the compact set $\mathcal{P}$, i.e. the full-distance
$$\kappa_N := \sup_{\mu\in\mathcal{P}}\mbox{dist}(\mu,\mathcal{P}_N)$$
tends to zero as $N$ to infinity, where $\mbox{dist}(\mu,\mathcal{P}_N) := \inf_{\mu' \in\mathcal{P}_N}\|\mu-\mu'\|_{\mathbb{R}^p}$. 

A crucial property for an efficiently computational procedure of reduced basis approaches is the parameter separability that can be defined as follows (cf.\ \cite{EPR10}). Such separability conditions, for instance, can be easily obtain by using the Empirical Interpolation Method (EIM) \cite{Bar04}, see also \cite{Ant2014}.
\begin{definition}
Assume that the functions $\sigma$, $\epsilon$, the desired control and state $\mb{u}_b$ and $\mb{E}_d$ admit the expansions
\begin{alignat*}{2}
 \sigma^{-1}(\cdot;\mu) &:= \sum_{q=1}^{Q^{\sigma}} \Theta_q^{\sigma}(\mu) {\sigma}^{-1}_q(\cdot),  &\quad \epsilon(\cdot;\mu) &:= \sum_{q=1}^{Q^{\epsilon}} \Theta_q^{\epsilon}(\mu) {\epsilon}_q(\cdot), \\	
 \mb{u}_d(\cdot;\mu) &:= \sum_{q=1}^{Q^{\mb{u}_d}} \Theta_q^{\mb{u}_d}(\mu) {\mb{u}_d}_q(\cdot), &\quad \mb{E}_d(\cdot;\mu) 
    &:= \sum_{q=1}^{Q^{\mb{E}_d}} \Theta_q^{\mb{E}_d}(\mu) {\mb{E}_d}_q(\cdot),
\end{alignat*}
where $Q^\sigma$, $Q^\epsilon$, $Q^{\mb{u}_d}$ and $Q^{\mb{E}_d}$ are finite positive integers, the functions $\Theta_q^{\sigma}, \Theta_q^{\epsilon}, \Theta_q^{\mb{u}_d}$ and $\Theta_q^{\mb{E}_d} : \mathcal{P} \to \mathbb{R}$, while the functions $\sigma^{-1}_q, \epsilon_q, {\mb{u}_d}_q : \Omega \to \mathbb{R}$ as well as $\mb{E}_d : D \to \mathbb{R}$ are independent of the parameter $\mu$.
\end{definition}
Due to \cref{11-2-18ct5} and \cref{8-5-19ct10}, without loss of generality we assume that
\begin{align*}
\sigma^{-1}_q \in [\overline{\sigma}^{-1}, \underline{\sigma}^{-1}], \quad 
\epsilon_q \in [\underline{\epsilon}, \overline{\epsilon}], \quad
\|{\mb{u}_d}_q\|_{\mb{L}^2(\Omega)} \le u_d
\quad\mbox{and}\quad
\|{\mb{E}_d}_q\|_{\mb{L}^2(D)} \le e_d
\end{align*}
for all the index $q$. 
We start with some auxiliary results.

\begin{lemma}\label{23-6-19ct1}
(i) For any given $\mu \in \mathcal{P}$ the inequality
\begin{align*}
\|S_{\mb{u}^1}(\mu)  -S_{\mb{u}^2}(\mu)\|_{\mathbf{H}(\mathbf{curl};\Omega)}  
\le C^\Omega_{\overline{\sigma}} \overline{\epsilon} \| u^1 - u^2\|_{\mb{L}^2(\Omega)}
\end{align*}
is satisfied for all $u^1, u^2 \in \mb{U}_{ad}$.

(ii) For any fixed $\mb{u} \in \mathbf{U}_{ad}$ the estimate
\begin{align*}
&\|S_{\mb{u}}(\mu^1)  -S_{\mb{u}}(\mu^2)\|_{\mathbf{H}(\mathbf{curl};\Omega)}  \\
&~\quad \le C^\Omega_{\overline{\sigma}} C_\mathbf{E} \underline{\sigma}^{-1}  \sum_{q=1}^{Q^{\sigma}}  |\Theta_q^{\sigma}(\mu^1) - \Theta_q^{\sigma}(\mu^2)|  + C^\Omega_{\overline{\sigma}}\overline{\epsilon} \|\overline{\mb{u}}\|_{\mathbb{R}^3} |\Omega|^{1/2} \sum_{q=1}^{Q^{\epsilon}} |\Theta_q^{\epsilon}(\mu^1) - \Theta_q^{\epsilon}(\mu^2)|
\end{align*}
holds true.
\end{lemma}

\begin{proof}
(i) For any $\mb{\Phi} \in \mathbf{H}_0(\mathbf{curl};\Omega)$ from the system \cref{11-2-18ct4} have that
\begin{align*}
\left(\sigma^{-1}(\mu)\nabla\times (S_{\mb{u}^1}(\mu)  -S_{\mb{u}^2}(\mu)), \nabla\times \mathbf{\Phi}\right)_{\mathbf{L}^2(\Omega)}
& = \left(\epsilon(\mu) (\mathbf{u}^1 - \mb{u}^2), \mathbf{\Phi}\right)_{\mathbf{L}^2(\Omega)}.
\end{align*}
Taking $\mb{\Phi} = S_{\mb{u}^1}(\mu)  -S_{\mb{u}^2}(\mu)$, with the aid of \cref{11-2-18ct5} and \cref{23-6-19ct2}, we obtain
\begin{align*}
\|S_{\mb{u}^1}(\mu)  -S_{\mb{u}^2}(\mu)\|^2_{\mathbf{H}(\mathbf{curl};\Omega)} \le C^\Omega_{\overline{\sigma}} \overline{\epsilon} \| u^1 - u^2\|_{\mb{L}^2(\Omega)} \| S_{\mb{u}^1}(\mu)  -S_{\mb{u}^2}(\mu)\|_{\mb{L}^2(\Omega)},
\end{align*}
which yields the desired inequality.

(ii) Likewise, we get
\begin{align*}
&\left(\sigma^{-1}(\mu^1)\nabla\times (S_{\mb{u}}(\mu^1)  -S_{\mb{u}}(\mu^2)), \nabla\times \mathbf{\Phi}\right)_{\mathbf{L}^2(\Omega)} \\
%&~\quad = \left(\sigma^{-1}(\mu^1)\nabla\times S_{\mb{u}}(\mu^1)  , \nabla\times \mathbf{\Phi}\right)_{\mathbf{L}^2(\Omega)} - \left(\sigma^{-1}(\mu^1)\nabla\times S_{\mb{u}}(\mu^2)  , \nabla\times \mathbf{\Phi}\right)_{\mathbf{L}^2(\Omega)}\\
&~\quad =\left(\epsilon(\mu^1) \mathbf{u}, \mathbf{\Phi}\right)_{\mathbf{L}^2(\Omega)} - \left(\sigma^{-1}(\mu^1)\nabla\times S_{\mb{u}}(\mu^2)  , \nabla\times \mathbf{\Phi}\right)_{\mathbf{L}^2(\Omega)}\\
%&~\quad =\left(\epsilon(\mu^2) \mathbf{u}, \mathbf{\Phi}\right)_{\mathbf{L}^2(\Omega)} + \left((\epsilon(\mu^1) - \epsilon(\mu^2)) \mathbf{u}, \mathbf{\Phi}\right)_{\mathbf{L}^2(\Omega)}\\
%&~\quad\quad - \left(\sigma^{-1}(\mu^1)\nabla\times S_{\mb{u}}(\mu^2)  , \nabla\times \mathbf{\Phi}\right)_{\mathbf{L}^2(\Omega)} \\
&~\quad = \left((\sigma^{-1}(\mu^2) - \sigma^{-1}(\mu^1))\nabla\times S_{\mb{u}}(\mu^2)  , \nabla\times \mathbf{\Phi}\right)_{\mathbf{L}^2(\Omega)} + \left((\epsilon(\mu^1) - \epsilon(\mu^2)) \mathbf{u}, \mathbf{\Phi}\right)_{\mathbf{L}^2(\Omega)}.
\end{align*}
Taking $\mb{\Phi} = S_{\mb{u}}(\mu^1)  -S_{\mb{u}}(\mu^2)$, we thus obtain
\begin{align*}
\|S_{\mb{u}}(\mu^1)  -S_{\mb{u}}(\mu^2)\|_{\mathbf{H}(\mathbf{curl};\Omega)}
& \le C^\Omega_{\overline{\sigma}}  \underline{\sigma}^{-1}\sum_{q=1}^{Q^{\sigma}} |\Theta_q^{\sigma}(\mu^1) - \Theta_q^{\sigma}(\mu^2)| \|S_{\mb{u}}(\mu^2)\|_{\mathbf{H}(\mathbf{curl};\Omega)} \\
&~\quad  + C^\Omega_{\overline{\sigma}}\overline{\epsilon}\sum_{q=1}^{Q^{\epsilon}} |\Theta_q^{\epsilon}(\mu^1) - \Theta_q^{\epsilon}(\mu^2)| \|\mb{u}\|_{\mathbf{L}^2(\Omega)}
\end{align*}
which together with \cref{3-5-19ct7} yield the desired inequality. The proof completes.
\end{proof}

\begin{lemma}\label{7-5-19ct1}
Let $\mb{u}^*_e(\mu)$ be the solution of the problem $(\mathbb{P}_e)$ with respect to the parameter $\mu\in\mathcal{P}$. Then, the estimate
\begin{align*}
&\|\mb{u}^*_e(\mu^1) - \mb{u}^*_e(\mu^2)\|_{\mb{L}^2(\Omega)}\\ 
&~\quad \le \sqrt{C^\sigma_{1/2}} \left( \sum_{q=1}^{Q^{\sigma}}  |\Theta_q^{\sigma}(\mu^1) - \Theta_q^{\sigma}(\mu^2)| \right)^{1/2} + \sqrt{C^\epsilon_{1/2}} \left( \sum_{q=1}^{Q^{\epsilon}} |\Theta_q^{\epsilon}(\mu^1) - \Theta_q^{\epsilon}(\mu^2)| \right)^{1/2} \\
&~\quad\quad + \sqrt{C^\epsilon_{1}} \sum_{q=1}^{Q^{\epsilon}} |\Theta_q^{\epsilon}(\mu^1) - \Theta_q^{\epsilon}(\mu^2)|  + \sqrt{C^{\mb{u}_d}_1} \left( \sum_{q=1}^{Q^{\mb{u}_d}} | \Theta_q^{\mb{u}_d}(\mu^1) - \Theta_q^{\mb{u}_d}(\mu^2)|^2\right)^{1/2}\\
&~\quad\quad + \sqrt{C^{\mb{E}_d}_{1/2}} \left( \sum_{q=1}^{Q^{\mb{E}_d}} | \Theta_q^{\mb{E}_d}(\mu^1) - \Theta_q^{\mb{E}_d}(\mu^2)|^2\right) ^{1/4}
\end{align*}
is satisfied for all $\mu^1, \mu^2 \in \mathcal{P}$, where
\begin{align*}
C^\sigma_{1/2} &:= 8C^\Omega_{\overline{\sigma}} C_\mathbf{E} \alpha^{-1} \underline{\sigma}^{-1}\underline{\epsilon}^{-1} \overline{\epsilon} \left( e_d + C_\mathbf{E}\right),  \\
C^\epsilon_{1/2} &:=  8\left( C^\Omega_{\overline{\sigma}}\overline{\epsilon} \left( e_d + C_\mathbf{E}\right) + C_\mathbf{F}\right) \alpha^{-1}\underline{\epsilon}^{-1}\overline{\epsilon}
\|\overline{\mb{u}}\|_{\mathbb{R}^3} |\Omega|^{1/2},\\
C^\epsilon_{1} &:=  4\left( (C^\Omega_{\overline{\sigma}})^2 \overline{\epsilon}^2 \left( e_d + C_\mathbf{E}\right)^2 + \left(\alpha \big(u_d + |\Omega|^{1/2}\|\overline{\mb{u}}\|_{\mathbb{R}^3} \big) + C_\mathbf{F}\right)^2\right) \alpha^{-2} \underline{\epsilon}^{-2}\overline{\epsilon}^2,\\
C^{\mb{u}_d}_1 &:= 4\underline{\epsilon}^{-2}\overline{\epsilon}^2 u^2_d Q^{\mb{u}_d} \quad \mbox{and}\\
C^{\mb{E}_d}_{1/2} &:= 8C_\mathbf{E}\alpha^{-1}  \underline{\epsilon}^{-1} \overline{\epsilon} e_d \left( Q^{\mb{E}_d}\right)^{1/2}.
\end{align*}
\end{lemma}

\begin{proof}
By the variational inequality \cref{24-4-19ct7}, we have
\begin{align*}
\left(\epsilon(\mu^1)\big(\mathbf{u}^*_e(\mu^2) - \mathbf{u}^*_e(\mu^1)\big), \mathbf{u}_d(\mu^1)  - \frac{1}{\alpha}\mathbf{F}^*_e(\mu^1) - \mathbf{u}^*_e(\mu^1)\right)_{\mathbf{L}^2(\Omega)} &\le 0 \\
\left(\epsilon(\mu^2)\big(\mathbf{u}^*_e(\mu^1) - \mathbf{u}^*_e(\mu^2)\big), \mathbf{u}_d(\mu^2)  - \frac{1}{\alpha}\mathbf{F}^*_e(\mu^2) - \mathbf{u}^*_e(\mu^2)\right)_{\mathbf{L}^2(\Omega)} &\le 0 
\end{align*}
which yield
\begin{align*}
&\alpha\underline{\epsilon}\|\mb{u}^*_e(\mu^2) - \mb{u}^*_e(\mu^1)\|^2_{\mb{L}^2(\Omega)} \\
&~\quad \le \left(\epsilon(\mu^1)\big(\mathbf{u}^*_e(\mu^2) - \mathbf{u}^*_e(\mu^1)\big), \alpha\big(\mb{u}_d(\mu^2) - \mb{u}_d(\mu^1)\big)\right)_{\mb{L}^2(\Omega)} \\
&~\qquad + \left(\big(\epsilon(\mu^2) - \epsilon(\mu^1)\big)\big(\mathbf{u}^*_e(\mu^2) - \mathbf{u}^*_e(\mu^1)\big), \alpha \mathbf{u}_d(\mu^2)  - \mathbf{F}^*_e(\mu^2) - \alpha \mathbf{u}^*_e(\mu^2)\right)_{\mb{L}^2(\Omega)}\\
&~\qquad + \left(\epsilon(\mu^1)\big(\mathbf{u}^*_e(\mu^2) - \mathbf{u}^*_e(\mu^1)\big), \mathbf{F}^*_e(\mu^1) - \mathbf{F}^*_e(\mu^2)\right)_{\mb{L}^2(\Omega)} 
:= I_1 + I_2 + I_3.
\end{align*}
We bound for the terms $I_1$, $I_2$ and $I_3$. First, we get
\begin{align*}
I_1
& \le \|\mb{u}^*_e(\mu^2) - \mb{u}^*_e(\mu^1)\|_{\mb{L}^2(\Omega)} \cdot \alpha\overline{\epsilon} u_d \left( Q^{\mb{u}_d}\right)^{1/2} \left( \sum_{q=1}^{Q^{\mb{u}_d}} \left| \Theta_q^{\mb{u}_d}(\mu^2) - \Theta_q^{\mb{u}_d}(\mu^1) \right|^2  \right)^{1/2} \\
& \le 4^{-1} \alpha \underline{\epsilon}\|\mb{u}^*_e(\mu^2) - \mb{u}^*_e(\mu^1)\|^2_{\mb{L}^2(\Omega)} + \alpha\underline{\epsilon}^{-1}\overline{\epsilon}^2 u^2_d Q^{\mb{u}_d} \sum_{q=1}^{Q^{\mb{u}_d}} \left| \Theta_q^{\mb{u}_d}(\mu^2) - \Theta_q^{\mb{u}_d}(\mu^1) \right|^2
\end{align*}
and, by \eqref{12-4-20ct1},		
\begin{align*}
I_2
& \le \left\|\big(\epsilon(\mu^2) - \epsilon(\mu^1)\big)\big(\mathbf{u}^*_e(\mu^2) - \mathbf{u}^*_e(\mu^1)\big)\right\|_{\mb{L}^2(\Omega)} \cdot\\
&~\quad \cdot\left(\alpha \|\mathbf{u}_d(\mu^2)\|_{\mb{L}^2(\Omega)}  + C_\mathbf{F} + \alpha \|\mathbf{u}^*_e(\mu^2)\|_{\mb{L}^2(\Omega)}\right)\\
&\le \overline{\epsilon}\sum_{q=1}^{Q^{\epsilon}} |\Theta_q^{\epsilon}(\mu^2) - \Theta_q^{\epsilon}(\mu^1)|
\|\mb{u}^*_e(\mu^2) - \mb{u}^*_e(\mu^1)\|_{\mb{L}^2(\Omega)}  \left(\alpha \big(u_d + |\Omega|^{1/2}\|\overline{\mb{u}}\|_{\mathbb{R}^3} \big) + C_\mathbf{F}\right)\\
&\le 4^{-1} \alpha \underline{\epsilon}\|\mb{u}^*_e(\mu^2) - \mb{u}^*_e(\mu^1)\|^2_{\mb{L}^2(\Omega)}\\
&~\quad + \alpha^{-1}\underline{\epsilon}^{-1}\overline{\epsilon}^2  \left(\alpha \big(u_d + |\Omega|^{1/2}\|\overline{\mb{u}}\|_{\mathbb{R}^3} \big) + C_\mathbf{F}\right)^2 \left(\sum_{q=1}^{Q^{\epsilon}} |\Theta_q^{\epsilon}(\mu^2) - \Theta_q^{\epsilon}(\mu^1)|\right)^2.
\end{align*}
For $I_3$ we write
\begin{align*}
I_3 
%&= \left(\epsilon(\mu^1) \mathbf{u}^*_e(\mu^1), \mathbf{F}^*_e(\mu^2) \right)_{\mb{L}^2(\Omega)} - \left(\epsilon(\mu^1)\mathbf{u}^*_e(\mu^2), \mathbf{F}^*_e(\mu^2) \right)_{\mb{L}^2(\Omega)} \\
%&~\quad + \left(\epsilon(\mu^1)\mathbf{u}^*_e(\mu^2), \mathbf{F}^*_e(\mu^1) \right)_{\mb{L}^2(\Omega)} - \left(\epsilon(\mu^1) \mathbf{u}^*_e(\mu^1), \mathbf{F}^*_e(\mu^1) \right)_{\mb{L}^2(\Omega)}\\
&= \left((\epsilon(\mu^1) - \epsilon(\mu^2))\mathbf{u}^*_e(\mu^1), \mathbf{F}^*_e(\mu^2) \right)_{\mb{L}^2(\Omega)} - \left((\epsilon(\mu^1)- \epsilon(\mu^2))\mathbf{u}^*_e(\mu^2), \mathbf{F}^*_e(\mu^2) \right)_{\mb{L}^2(\Omega)}\\
&~\quad + \left(\epsilon(\mu^2) \mathbf{u}^*_e(\mu^1), \mathbf{F}^*_e(\mu^2) \right)_{\mb{L}^2(\Omega)} - \left(\epsilon(\mu^2)\mathbf{u}^*_e(\mu^2), \mathbf{F}^*_e(\mu^2) \right)_{\mb{L}^2(\Omega)}\\
&~\quad + \left(\epsilon(\mu^1) (\mathbf{E}^*_e(\mu^1) - \mathbf{E}_d(\mu^1)), S_{\mb{u}^*_e(\mu^2)}(\mu^1) - \mathbf{E}^*_e(\mu^1)\right)_{\mathbf{L}^2(D)}\\
&\le 2\overline{\epsilon}C_\mathbf{F}|\Omega|^{1/2}\|\overline{u}\|_{\mathbb{R}^3}\sum_{q=1}^{Q^{\epsilon}} |\Theta_q^{\epsilon}(\mu^2) - \Theta_q^{\epsilon}(\mu^1)|\\
&~\quad + \left(\epsilon(\mu^2) (\mathbf{E}^*_e(\mu^2) - \mathbf{E}_d(\mu^2)), S_{\mb{u}^*_e(\mu^1)}(\mu^2)- \mathbf{E}^*_e(\mu^2)\right)_{\mathbf{L}^2(D)}\\
&~\quad + \left(\epsilon(\mu^1) (\mathbf{E}^*_e(\mu^1) - \mathbf{E}_d(\mu^1)), S_{\mb{u}^*_e(\mu^2)}(\mu^1) - \mathbf{E}^*_e(\mu^1)\right)_{\mathbf{L}^2(D)} := I_3^1+I_3^2+I_3^3
\end{align*}
with
\begin{align*}
I_3^2 + I_3^3
%&= \left((\epsilon(\mu^1) - \epsilon(\mu^2) (\mathbf{E}^*_e(\mu^1) - \mathbf{E}_d(\mu^1)), S_{\mb{u}^*_e(\mu^2)}(\mu^1) - \mathbf{E}^*_e(\mu^1)\right)_{\mathbf{L}^2(D)} \notag\\
%&~\quad +  \left(\epsilon(\mu^2) (\mathbf{E}^*_e(\mu^1) - \mathbf{E}_d(\mu^1)), S_{\mb{u}^*_e(\mu^2)}(\mu^1) - \mathbf{E}^*_e(\mu^1)\right)_{\mathbf{L}^2(D)}\notag\\
%&~\quad +  \left(\epsilon(\mu^2) (\mathbf{E}^*_e(\mu^2) - \mathbf{E}_d(\mu^2)), S_{\mb{u}^*_e(\mu^1)}(\mu^2)- \mathbf{E}^*_e(\mu^2) \right)_{\mathbf{L}^2(D)}\\
&= \left((\epsilon(\mu^1) - \epsilon(\mu^2) (\mathbf{E}^*_e(\mu^1) - \mathbf{E}_d(\mu^1)), S_{\mb{u}^*_e(\mu^2)}(\mu^1) - \mathbf{E}^*_e(\mu^1)\right)_{\mathbf{L}^2(D)} \notag\\
&~\quad + \left(\epsilon(\mu^2) (\mathbf{E}^*_e(\mu^1) - \mathbf{E}_d(\mu^1)), S_{\mb{u}^*_e(\mu^2)}(\mu^1) - S_{\mb{u}^*_e(\mu^2)}(\mu^2)\right)_{\mathbf{L}^2(D)} \notag\\
&~\quad  +  \left(\epsilon(\mu^2) (\mathbf{E}^*_e(\mu^2) - \mathbf{E}_d(\mu^2)), S_{\mb{u}^*_e(\mu^1)}(\mu^2)- S_{\mb{u}^*_e(\mu^1)}(\mu^1) \right)_{\mathbf{L}^2(D)}\notag\\
&~\quad  +  \left(\epsilon(\mu^2) (\mathbf{E}^*_e(\mu^1) - \mathbf{E}_d(\mu^1)), S_{\mb{u}^*_e(\mu^2)}(\mu^2) - \mb{E}^*_e(\mu^1)\right)_{\mathbf{L}^2(D)} \notag\\
&~\quad  + \left(\epsilon(\mu^2) (\mathbf{E}^*_e(\mu^2) - \mathbf{E}_d(\mu^2)), S_{\mb{u}^*_e(\mu^1)}(\mu^1) - \mb{E}^*_e(\mu^2) \right)_{\mathbf{L}^2(D)} \\
& := J_1+J_2+J_3+J_4+J_5.
\end{align*}
Since $\mb{E}^*_e(\mu^1) = S_{\mb{u}^*_e(\mu^1)}(\mu^1)$,  we with the aid of \cref{8-5-19ct10}, \cref{3-5-19ct7} and \cref{23-6-19ct1} get
\begin{align*}
J_1 
&\le \overline{\epsilon} \sum_{q=1}^{Q^{\epsilon}} |\Theta_q^{\epsilon}(\mu^1) - \Theta_q^{\epsilon}(\mu^2)| \cdot\notag\\
&~\quad\quad\cdot \left( \|\mathbf{E}^*_e(\mu^1) \|_{\mathbf{L}^2(D)}  + \|\mathbf{E}_d(\mu^1)\|_{\mb{L}^2(D)} \right) \cdot \left( \|S_{\mb{u}^*_e(\mu^2)}(\mu^1) - S_{\mb{u}^*_e(\mu^1)}(\mu^1)\|_{\mathbf{L}^2(D)}\right) \notag\\
&\le C^\Omega_{\overline{\sigma}} \overline{\epsilon}^2\left( e_d + C_\mathbf{E}\right)   \sum_{q=1}^{Q^{\epsilon}} |\Theta_q^{\epsilon}(\mu^1) - \Theta_q^{\epsilon}(\mu^2)| \|\mb{u}^*_e(\mu^2) - \mb{u}^*_e(\mu^1)\|_{\mb{L}^2(\Omega)} \notag\\
&\le 4^{-1}\alpha \underline{\epsilon} \|\mb{u}^*_e(\mu^2) - \mb{u}^*_e(\mu^1)\|^2_{\mb{L}^2(\Omega)} \notag\\
&~\quad + (C^\Omega_{\overline{\sigma}})^2\alpha^{-1} \underline{\epsilon}^{-1} \overline{\epsilon}^4 \left( e_d + C_\mathbf{E}\right)^2 \left(\sum_{q=1}^{Q^{\epsilon}} |\Theta_q^{\epsilon}(\mu^1) - \Theta_q^{\epsilon}(\mu^2)|\right)^2
\end{align*}
and
\begin{align*}
J_2+J_3 
%&\le \overline{\epsilon} \left( e_d + C_\mathbf{E}\right) \cdot \notag\\
%&~\quad\quad\cdot \left( \|S_{\mb{u}^*_e(\mu^2)}(\mu^1) - S_{\mb{u}^*_e(\mu^2)}(\mu^2)\|_{\mb{L}^2(\Omega)} + \|S_{\mb{u}^*_e(\mu^1)}(\mu^1) - S_{\mb{u}^*_e(\mu^1)}(\mu^2)\|_{\mb{L}^2(\Omega)}\right)\\
&\le 2C^\Omega_{\overline{\sigma}} C_\mathbf{E} \underline{\sigma}^{-1} \overline{\epsilon} \left( e_d + C_\mathbf{E}\right) \sum_{q=1}^{Q^{\sigma}}  |\Theta_q^{\sigma}(\mu^1) - \Theta_q^{\sigma}(\mu^2)|  \\
&~\quad + 2C^\Omega_{\overline{\sigma}}\overline{\epsilon}^2 \left( e_d + C_\mathbf{E}\right) \|\overline{\mb{u}}\|_{\mathbb{R}^3} |\Omega|^{1/2} \sum_{q=1}^{Q^{\epsilon}} |\Theta_q^{\epsilon}(\mu^1) - \Theta_q^{\epsilon}(\mu^2)|
\end{align*}
and
\begin{align*}
J_4 +J_5  
&=  \left(\epsilon(\mu^2) (\mathbf{E}^*_e(\mu^1) - \mathbf{E}_d(\mu^1)), \mb{E}^*_e(\mu^2) - \mb{E}^*_e(\mu^1)\right)_{\mathbf{L}^2(D)} \notag\\
&~\quad  +  \left(\epsilon(\mu^2) (\mathbf{E}^*_e(\mu^2) - \mathbf{E}_d(\mu^1)), \mb{E}^*_e(\mu^1) - \mb{E}^*_e(\mu^2) \right)_{\mathbf{L}^2(D)} \notag\\
&~\quad  + \left(\epsilon(\mu^2) (\mathbf{E}_d(\mu^1) - \mathbf{E}_d(\mu^2)), \mb{E}^*_e(\mu^1) - \mb{E}^*_e(\mu^2) \right)_{\mathbf{L}^2(D)} \notag\\
& \le \left(\epsilon(\mu^2) (\mathbf{E}_d(\mu^1) - \mathbf{E}_d(\mu^2)), \mb{E}^*_e(\mu^1) - \mb{E}^*_e(\mu^2) \right)_{\mathbf{L}^2(D)} \notag\\
&~\quad - \underline{\epsilon} \|\mathbf{E}^*_e(\mu^2) - \mathbf{E}^*_e(\mu^1) \|^2_{\mathbf{L}^2(D)} \notag\\
& \le \overline{\epsilon}  \left( \|\mathbf{E}^*_e(\mu^1) \|_{\mathbf{L}^2(D)} + \mathbf{E}^*_e(\mu^2) \|_{\mathbf{L}^2(D)}\right) \|\mathbf{E}_d(\mu^1) - \mathbf{E}_d(\mu^2)\|_{\mb{L}^2(D)} \notag\\
& \le 2 C_\mathbf{E} \overline{\epsilon} e_d\left( Q^{\mb{E}_d}\right)^{1/2} \Big( \sum_{q=1}^{Q^{\mb{E}_d}} | \Theta_q^{\mb{E}_d}(\mu^1) - \Theta_q^{\mb{E}_d}(\mu^2)|^2 \Big)^{1/2}.
\end{align*}
Therefore, we arrive at
\begin{align*}
I_3 
&\le 2\overline{\epsilon}C_\mathbf{F}|\Omega|^{1/2}\|\overline{u}\|_{\mathbb{R}^3}\sum_{q=1}^{Q^{\epsilon}} |\Theta_q^{\epsilon}(\mu^2) - \Theta_q^{\epsilon}(\mu^1)|\\
&~\quad + 4^{-1}\alpha \underline{\epsilon} \|\mb{u}^*_e(\mu^2) - \mb{u}^*_e(\mu^1)\|^2_{\mb{L}^2(\Omega)} \notag\\
&~\quad + (C^\Omega_{\overline{\sigma}})^2\alpha^{-1} \underline{\epsilon}^{-1} \overline{\epsilon}^4 \left( e_d + C_\mathbf{E}\right)^2 \left(\sum_{q=1}^{Q^{\epsilon}} |\Theta_q^{\epsilon}(\mu^1) - \Theta_q^{\epsilon}(\mu^2)|\right)^2 \\
&~\quad + 2C^\Omega_{\overline{\sigma}} C_\mathbf{E} \underline{\sigma}^{-1} \overline{\epsilon} \left( e_d + C_\mathbf{E}\right) \sum_{q=1}^{Q^{\sigma}}  |\Theta_q^{\sigma}(\mu^1) - \Theta_q^{\sigma}(\mu^2)|  \\
&~\quad + 2C^\Omega_{\overline{\sigma}}\overline{\epsilon}^2 \left( e_d + C_\mathbf{E}\right) \|\overline{\mb{u}}\|_{\mathbb{R}^3} |\Omega|^{1/2} \sum_{q=1}^{Q^{\epsilon}} |\Theta_q^{\epsilon}(\mu^1) - \Theta_q^{\epsilon}(\mu^2)|\\
&~\quad + 2 C_\mathbf{E} \overline{\epsilon} e_d\left( Q^{\mb{E}_d}\right)^{1/2} \Big( \sum_{q=1}^{Q^{\mb{E}_d}} | \Theta_q^{\mb{E}_d}(\mu^1) - \Theta_q^{\mb{E}_d}(\mu^2)|^2 \Big)^{1/2}.
\end{align*}
The desired estimate follows from the bounds for $I_1,~ I_2,~ I_3$ and \cref{23-6-19ct1}, which finishes the proof.
\end{proof}

Now we state the similar results for the finite dimensional approximation problem $(\mathbb{P}_h)$ and the reduced basis approach $(\mathbb{P}_N)$, their proofs follow exactly as in the continuous case $(\mathbb{P}_e)$, therefore omitted here.

\begin{lemma}\label{8-5-19ct12}
Let $\mb{u}^*_h(\mu)$ and $\mb{u}^*_N(\mu)$ respectively be the solution of the problems $(\mathbb{P}_h)$ and $(\mathbb{P}_N)$ at the given parameter $\mu\in\mathcal{P}$. Then, the estimates
\begin{align*}
&\|\mb{u}^*_h(\mu^1) - \mb{u}^*_h(\mu^2)\|_{\mb{L}^2(\Omega)}\\ 
&~\quad \le \sqrt{C^\sigma_{1/2}} \left( \sum_{q=1}^{Q^{\sigma}}  |\Theta_q^{\sigma}(\mu^1) - \Theta_q^{\sigma}(\mu^2)| \right)^{1/2} + \sqrt{C^\epsilon_{1/2}} \left( \sum_{q=1}^{Q^{\epsilon}} |\Theta_q^{\epsilon}(\mu^1) - \Theta_q^{\epsilon}(\mu^2)| \right)^{1/2} \\
&~\quad\quad + \sqrt{C^\epsilon_{1}} \sum_{q=1}^{Q^{\epsilon}} |\Theta_q^{\epsilon}(\mu^1) - \Theta_q^{\epsilon}(\mu^2)|  + \sqrt{C^{\mb{u}_d}_1} \left( \sum_{q=1}^{Q^{\mb{u}_d}} | \Theta_q^{\mb{u}_d}(\mu^1) - \Theta_q^{\mb{u}_d}(\mu^2)|^2\right)^{1/2}\\
&~\quad\quad + \sqrt{C^{\mb{E}_d}_{1/2}} \left( \sum_{q=1}^{Q^{\mb{E}_d}} | \Theta_q^{\mb{E}_d}(\mu^1) - \Theta_q^{\mb{E}_d}(\mu^2)|^2\right) ^{1/4}
\end{align*}
and
\begin{align*}
&\|\mb{u}^*_N(\mu^1) - \mb{u}^*_N(\mu^2)\|_{\mb{L}^2(\Omega)}\\ 
&~\quad \le \sqrt{C^\sigma_{1/2}} \left( \sum_{q=1}^{Q^{\sigma}}  |\Theta_q^{\sigma}(\mu^1) - \Theta_q^{\sigma}(\mu^2)| \right)^{1/2} + \sqrt{C^\epsilon_{1/2}} \left( \sum_{q=1}^{Q^{\epsilon}} |\Theta_q^{\epsilon}(\mu^1) - \Theta_q^{\epsilon}(\mu^2)| \right)^{1/2} \\
&~\quad\quad + \sqrt{C^\epsilon_{1}} \sum_{q=1}^{Q^{\epsilon}} |\Theta_q^{\epsilon}(\mu^1) - \Theta_q^{\epsilon}(\mu^2)|  + \sqrt{C^{\mb{u}_d}_1} \left( \sum_{q=1}^{Q^{\mb{u}_d}} | \Theta_q^{\mb{u}_d}(\mu^1) - \Theta_q^{\mb{u}_d}(\mu^2)|^2\right)^{1/2}\\
&~\quad\quad + \sqrt{C^{\mb{E}_d}_{1/2}} \left( \sum_{q=1}^{Q^{\mb{E}_d}} | \Theta_q^{\mb{E}_d}(\mu^1) - \Theta_q^{\mb{E}_d}(\mu^2)|^2\right) ^{1/4}
\end{align*}
hold true for all $\mu^1, \mu^2 \in \mathcal{P}$.
\end{lemma}

We are in the position to state the main result of this section on the uniform convergence of reduced order solutions. To do so, we assume $\mb{u}^*_h(\mu) = \mb{u}^*_N(\mu)$ for $\mu$ belonging to the parameter sample $\mathcal{P}_N$. For the state equation, this assumption is  the basic consistency property of an Reduced Basis scheme, which simply put is the reproduction of solutions (cf.\ \cite[Proposition~2.20]{Ha17}). For a justification of this assumption for optimal control problems, see \cite[pp.~A283]{Ali}.

\begin{theorem}\label{8-5-19ct17}
Assume that the functions $\Theta_q^{\sigma}, \Theta_q^{\epsilon}, \Theta_q^{\mb{u}_d}, \Theta_q^{\mb{E}_d} : \mathcal{P} \to \mathbb{R}$ are H\"older continuous, i.e. 
\begin{align*}
|\Theta_q^{\sigma}(\mu^1) - \Theta_q^{\sigma}(\mu^2)| 
&\le L^\sigma \|\mu^1 - \mu^2\|^{\gamma^\sigma}_{\mathbb{R}^p} \\
|\Theta_q^{\epsilon}(\mu^1) - \Theta_q^{\epsilon}(\mu^2)| 
&\le L^\epsilon \|\mu^1 - \mu^2\|^{\gamma^\epsilon}_{\mathbb{R}^p}\\
|\Theta_q^{\mb{u}_d}(\mu^1) - \Theta_q^{\mb{u}_d}(\mu^2)| 
&\le L^{\mb{u}_d} \|\mu^1 - \mu^2\|^{\gamma^{\mb{u}_d}}_{\mathbb{R}^p}\\
|\Theta_q^{\mb{E}_d}(\mu^1) - \Theta_q^{\mb{E}_d}(\mu^2)| 
&\le L^{\mb{E}_d} \|\mu^1 - \mu^2\|^{\gamma^{\mb{E}_d}}_{\mathbb{R}^p} 
\end{align*}
for all $\mu^1, \mu^2 \in\mathcal{P}$ and all the index $q$ with some positive constants $L^\sigma$, $L^\epsilon$, $L^{\mb{u}_d}$, $L^{\mb{E}_d}$ and $\gamma^\sigma$, $\gamma^\epsilon$, $\gamma^{\mb{u}_d}$, $\gamma^{\mb{E}_d}$. For any given $\mu\in\mathcal{P}$ let $\mb{u}^*_h(\mu)$ and $\mb{u}^*_N(\mu)$ be the solutions of the problems $(\mathbb{P}_h)$ and $(\mathbb{P}_N)$, respectively. Then the estimate
$$\|\mb{u}^*_h(\mu) - \mb{u}^*_N(\mu)\|_{\mb{L}^2(\Omega)} \le C\kappa^\gamma_N$$
is established, where $\gamma := \frac{1}{2}\min \left( \gamma^\sigma, \gamma^\epsilon, 2\gamma^{\mb{u}_d}, \gamma^{\mb{E}_d}\right)  >0$.
\end{theorem}

\begin{proof}
For all $\mu^1, \mu^2 \in \mathcal{P}$ we deduce from  \Cref{8-5-19ct12} that
\begin{align*}
&\max\left( \|\mb{u}^*_h(\mu^1) - \mb{u}^*_h(\mu^2)\|_{\mb{L}^2(\Omega)}, \|\mb{u}^*_N(\mu^1) - \mb{u}^*_N(\mu^2)\|^2_{\mb{L}(\Omega)}\right) \\
&~\quad \le C\left( \|\mu^1 - \mu^2\|^{\gamma^\sigma/2}_{\mathbb{R}^p} + \|\mu^1 - \mu^2\|^{\gamma^\epsilon/2}_{\mathbb{R}^p} + \|\mu^1 - \mu^2\|^{\gamma^{\mb{u}_d}} + \|\mu^1 - \mu^2\|^{\gamma^{\mb{E}_d}/2}_{\mathbb{R}^p}\right),
\end{align*}
where the positive constant $C$ is independent of the parameters. For any fixed $\mu \in \mathcal{P}$, since the set $\mathcal{P}_N$ is finite, there exists $\mu^* \in \arg \min_{\mu' \in \mathcal{P}_N}\|\mu-\mu'\|_{\mathbb{R}^p}$. By $\mu^* \in \mathcal{P}_N$, we get $\mb{u}^*_h(\mu^*) = \mb{u}^*_N(\mu^*)$  and therefore obtain that
\begin{align*}
&\|\mb{u}^*_h(\mu) - \mb{u}^*_N(\mu)\|_{\mb{L}^2(\Omega)} \\
&~\quad = \|\mb{u}^*_h(\mu) - \mb{u}^*_h(\mu^*)  + \mb{u}^*_N(\mu^*) - \mb{u}^*_N(\mu)\|_{\mb{L}^2(\Omega)} \\
&~\quad \le  \|\mb{u}^*_h(\mu) - \mb{u}^*_h(\mu^*) \|_{\mb{L}^2(\Omega)} + \|\mb{u}^*_N(\mu^*) - \mb{u}^*_N(\mu)\|_{\mb{L}^2(\Omega)}\\
&~\quad \le C\left( \|\mu - \mu^*\|^{\gamma^\sigma/2}_{\mathbb{R}^p} + \|\mu - \mu^*\|^{\gamma^\epsilon/2}_{\mathbb{R}^p} + \|\mu^1 - \mu^2\|^{\gamma^{\mb{u}_d}} + \|\mu^1 - \mu^2\|^{\gamma^{\mb{E}_d}/2}_{\mathbb{R}^p}\right) \\
&~\quad \le C\kappa^\gamma_N,
\end{align*}
which finishes the proof.
\end{proof}

\section{A posteriori error estimation} \label{Post-estimate}

In the greedy sampling procedure, a possibility for the error estimator is that 
$$\Delta_N(\mathcal{E}_N,V_N;\mu) = \|\mathbf{u}^*_h(\mu) - \mathbf{u}^*_N(\mu)\|_{\mathbf{L}^2(\Omega)}.$$ 
However, this estimator depends on $\mathbf{u}^*_h(\mu)$, i.e. the high dimensional problem $(\mathbb{P}_h)$. In view of a posteriori error estimates we wish to construct an error estimator which is independent of the solution to $(\mathbb{P}_h)$.

For a given $\mu\in \mathcal{P}$ let $(\mb{u}^*_N(\mu),\mb{E}^*_N(\mu),\mb{F}^*_N(\mu))$ satisfy the system \cref{24-4-19ct3**}--\cref{24-4-19ct7**}. Assume that $\widehat{\mb{E}}_h (\mu) \in \mathcal{E}_h$ and $\widehat{\mb{F}}_h (\mu) \in \mathcal{E}_h$ are defined by
\begin{equation}\label{29-4-19ct1}
\begin{cases}
\left(\sigma^{-1}(\mu)\nabla\times \widehat{\mb{E}}_h (\mu), \nabla\times \mathbf{\Phi}_h\right)_{\mathbf{L}^2(\Omega)} &= \left(\epsilon(\mu) \mathbf{u}^*_N(\mu), \mathbf{\Phi}_h\right)_{\mathbf{L}^2(\Omega)} \\
\left(\epsilon(\mu) \widehat{\mb{E}}_h (\mu), \nabla\phi_h\right)_{\mathbf{L}^2(\Omega)} &= -(\rho(\mu),\phi_h)_{L^2(\Omega)} 
\end{cases}
\end{equation}
and
\begin{equation}\label{29-4-19ct2}
\begin{cases}
\left(\sigma^{-1}(\mu)\nabla\times \widehat{\mb{F}}_h(\mu), \nabla\times \mathbf{\Phi}_h\right)_{\mathbf{L}^2(\Omega)} &= \left(\epsilon(\mu) (\mathbf{E}^*_N(\mu) - \mathbf{E}_d(\mu)), \mathbf{\Phi}_h\right)_{\mathbf{L}^2(D)}\\
\left(\epsilon(\mu) \widehat{\mb{F}}_h(\mu), \nabla\phi_h\right)_{\mathbf{L}^2(\Omega)} &= 0
\end{cases}
\end{equation}
for all $\left( \mb{\Phi}_h, \phi_h\right) \in \mathcal{E}_h \times V_h$, respectively. We further define the residuals $R_{\mb{E}} := R_{\mb{E}}(\cdot;\mu) \in \mathcal{E}_h^*$ and $R_{\mb{F}} := R_{\mb{F}}(\cdot;\mu) \in \mathcal{E}_h^*$ via the following identities
\begin{align*}
&R_{\mb{E}}(\mathbf{\Phi}_h;\mu) := \left(\epsilon(\mu) \mathbf{u}^*_N(\mu), \mathbf{\Phi}_h\right)_{\mathbf{L}^2(\Omega)} - \left(\sigma^{-1}(\mu)\nabla\times \mb{E}^*_N(\mu), \nabla\times \mathbf{\Phi}_h\right)_{\mathbf{L}^2(\Omega)} \\
&R_{\mb{F}}(\mathbf{\Phi}_h;\mu) := \left(\epsilon(\mu) (\mathbf{E}^*_N(\mu) - \mathbf{E}_d(\mu)), \mathbf{\Phi}_h\right)_{\mathbf{L}^2(D)} - \left(\sigma^{-1}(\mu)\nabla\times \mb{F}^*_N(\mu), \nabla\times \mathbf{\Phi}_h\right)_{\mathbf{L}^2(\Omega)}  
\end{align*}
for all $\mathbf{\Phi}_h\in\mathcal{E}_h$. 
%We note that
%$$R_{\mb{E}}(\cdot;\mu)_{|\mathcal{E}_N} \equiv 0 \quad\mbox{and}\quad R_{\mb{F}}(\cdot;\mu)_{|\mathcal{E}_N} \equiv 0$$
%due to the equations \cref{24-4-19ct3**} and \cref{24-4-19ct5**}, respectively.

To begin, we present some auxiliary results.

\begin{lemma}\label{30-4-19}
Let $(\mb{u}^*_h(\mu),\mb{E}^*_h(\mu),\mb{F}^*_h(\mu))$ and $(\mb{u}^*_N(\mu),\mb{E}^*_N(\mu),\mb{F}^*_N(\mu))$ respectively satisfy the systems \cref{24-4-19ct3*}--\cref{24-4-19ct7*} and \cref{24-4-19ct3**}--\cref{24-4-19ct7**}  at a given $\mu\in \mathcal{P}$. Then the following inequalities
\begin{align*}
\|\mb{E}^*_h(\mu) - \widehat{\mb{E}}_h(\mu)\|_{\mathbf{H}(\mathbf{curl};\Omega)}
&\le C^\Omega_{\overline{\sigma}}  \overline{\epsilon}\|\mb{u}^*_h(\mu) - \mathbf{u}^*_N(\mu)\|_{\mathbf{L}^2(\Omega)} \\
 \|\mb{F}^*_h(\mu) - \widehat{\mb{F}}_h(\mu)\|_{\mathbf{H}(\mathbf{curl};\Omega)}
&\le C^\Omega_{\overline{\sigma}}  \overline{\epsilon}\|\mb{E}^*_h(\mu) - \mathbf{E}^*_N(\mu) \|_{\mathbf{L}^2(\Omega)} 
\end{align*}
are satisfied.
\end{lemma}

\begin{proof}
By \cref{30-4-19ct1}, we have
\begin{align*}
\|\mb{E}^*_h - \widehat{\mb{E}}_h\|^2_{\mathbf{H}(\mathbf{curl};\Omega)}
&\le C^\Omega_{\overline{\sigma}} (\sigma^{-1}\nabla\times (\mb{E}^*_h - \widehat{\mb{E}}_h), \nabla\times (\mb{E}^*_h - \widehat{\mb{E}}_h))_{\mathbf{L}^2(\Omega)}\\
&= C^\Omega_{\overline{\sigma}} \left(\epsilon (\mb{u}^*_h - \mathbf{u}^*_N), \mb{E}^*_h - \widehat{\mb{E}}_h\right)_{\mathbf{L}^2(\Omega)}\\
&\le C^\Omega_{\overline{\sigma}}  \overline{\epsilon}\|\mb{u}^*_h - \mathbf{u}^*_N\|_{\mathbf{L}^2(\Omega)}  \|\mb{E}^*_h - \widehat{\mb{E}}_h\|_{\mathbf{H}(\mathbf{curl};\Omega)}
\end{align*}
which implies the first inequality. Likewise, we get
\begin{align*}
\|\mb{F}^*_h - \widehat{\mb{F}}_h\|^2_{\mathbf{H}(\mathbf{curl};\Omega)}
&\le C^\Omega_{\overline{\sigma}} (\sigma^{-1}\nabla\times (\mb{F}^*_h - \widehat{\mb{F}}_h), \nabla\times (\mb{F}^*_h - \widehat{\mb{F}}_h))_{\mathbf{L}^2(\Omega)}\\
&= C^\Omega_{\overline{\sigma}} \left(\epsilon (\mb{E}^*_h - \mathbf{E}^*_N), \mb{F}^*_h - \widehat{\mb{F}}_h\right)_{\mathbf{L}^2(D)}\\
&\le C^\Omega_{\overline{\sigma}} \overline{\epsilon}\|\mb{E}^*_h - \mathbf{E}^*_N \|_{\mathbf{L}^2(\Omega)}  \|\mb{F}^*_h - \widehat{\mb{F}}_h\|_{\mathbf{H}(\mathbf{curl};\Omega)}.
\end{align*}
The proof completes.
\end{proof}

\begin{lemma}\label{2-5-19}
The inequalities
\begin{align*}
\underline{\sigma} \|R_{\mb{E}}(\cdot;\mu)\|_{{_{\mathbf{H}(\mathbf{curl};\Omega)}}^*} 
&\le \|\mb{E}^*_N(\mu) - \widehat{\mb{E}}_h(\mu)\|_{\mathbf{H}(\mathbf{curl};\Omega)} \le C^\Omega_{\overline{\sigma}} \|R_{\mb{E}}(\cdot;\mu)\|_{{_{\mathbf{H}(\mathbf{curl};\Omega)}}^*} \\
\underline{\sigma} \|R_{\mb{F}}(\cdot;\mu)\|_{{_{\mathbf{H}(\mathbf{curl};\Omega)}}^*} 
&\le \|\mb{F}^*_N(\mu) - \widehat{\mb{F}}_h(\mu)\|_{\mathbf{H}(\mathbf{curl};\Omega)} \le C^\Omega_{\overline{\sigma}}\|R_{\mb{F}}(\cdot;\mu)\|_{{_{\mathbf{H}(\mathbf{curl};\Omega)}}^*}
\end{align*}
hold true.
\end{lemma}

\begin{proof}
We have
\begin{align*}
 \|\mb{E}^*_N - \widehat{\mb{E}}_h\|^2_{\mathbf{H}(\mathbf{curl};\Omega)}& \le C^\Omega_{\overline{\sigma}} \left(\sigma^{-1}\nabla\times (\mb{E}^*_N - \widehat{\mb{E}}_h), \nabla\times (\mb{E}^*_N - \widehat{\mb{E}}_h) \right)_{\mathbf{L}^2(\Omega)}\notag\\
%& = C^\Omega_{\overline{\sigma}} \left(\sigma^{-1}\nabla\times \widehat{\mb{E}}_h, \nabla\times (\widehat{\mb{E}}_h - \mb{E}^*_N) \right)_{\mathbf{L}^2(\Omega)} \notag\\
%&~\quad - C^\Omega_{\overline{\sigma}}\left(\sigma^{-1}\nabla\times \mb{E}^*_N, \nabla\times (\widehat{\mb{E}}_h - \mb{E}^*_N) \right)_{\mathbf{L}^2(\Omega)}\\
& = C^\Omega_{\overline{\sigma}}\left(\epsilon \mathbf{u}^*_N, \widehat{\mb{E}}_h - \mb{E}^*_N\right)_{\mathbf{L}^2(\Omega)}\\
&~\quad - C^\Omega_{\overline{\sigma}}\left(\sigma^{-1}\nabla\times \mb{E}^*_N, \nabla\times (\widehat{\mb{E}}_h - \mb{E}^*_N) \right)_{\mathbf{L}^2(\Omega)} \notag\\
&= C^\Omega_{\overline{\sigma}} R_{\mb{E}}(\widehat{\mb{E}}_h - \mb{E}^*_N)\\
& \le C^\Omega_{\overline{\sigma}}  \|R_{\mb{E}}\|_{{_{\mathbf{H}(\mathbf{curl};\Omega)}}^*} \|\widehat{\mb{E}}_h - \mb{E}^*_N\|_{\mathbf{H}(\mathbf{curl};\Omega)}.
\end{align*}
To show the lower bound, we first take $r_{\mb{E}} \in \mathcal{E}_h$ such that
$$R_{\mb{E}}(\mb{\Phi}_h) = (r_{\mb{E}};\mb{\Phi}_h)_{\mathbf{H}(\mathbf{curl};\Omega)},~ \forall \mb{\Phi}_h \in\mathcal{E}_h \mbox{~and~} \|r_{\mb{E}}\|_{\mathbf{H}(\mathbf{curl};\Omega)} = \|R_{\mb{E}}\|_{{_{\mathbf{H}(\mathbf{curl};\Omega)}}^*},$$
by Riesz Representation Theorem. In view of the above argument, we arrive at 
\begin{align*}
\|R_{\mb{E}}\|^2_{{\mathbf{H}(\mathbf{curl};\Omega)}^*}
&= (r_{\mb{E}},r_{\mb{E}})_{\mathbf{H}(\mathbf{curl};\Omega)}\\
&= R_{\mb{E}}(r_{\mb{E}})\\
&= \left(\epsilon \mathbf{u}^*_N, r_{\mb{E}}\right)_{\mathbf{L}^2(\Omega)} - \left(\sigma^{-1}\nabla\times \mb{E}^*_N, \nabla\times r_{\mb{E}}\right)_{\mathbf{L}^2(\Omega)}\\
&= \left(\sigma^{-1}\nabla\times (\widehat{\mb{E}}_h - \mb{E}^*_N), \nabla\times r_{\mb{E}}\right)_{\mathbf{L}^2(\Omega)}\\ 
&\le \underline{\sigma}^{-1}\|r_{\mb{E}}\|_{\mathbf{H}(\mathbf{curl};\Omega)} \|\widehat{\mb{E}}_h - \mb{E}^*_N\|_{\mathbf{H}(\mathbf{curl};\Omega)}\\
&= \underline{\sigma}^{-1}\|R_{\mb{E}}\|_{{_{\mathbf{H}(\mathbf{curl};\Omega)}}^*} \|\widehat{\mb{E}}_h - \mb{E}^*_N\|_{\mathbf{H}(\mathbf{curl};\Omega)}.
\end{align*}
The remaining inequalities follows by the same arguments, therefore omitted here.
\end{proof}

We now state the main results of the section.

\begin{theorem}\label{2-5-19-DL}
Let $(\mb{u}^*_h(\mu),\mb{E}^*_h(\mu),\mb{F}^*_h(\mu))$ and $(\mb{u}^*_N(\mu),\mb{E}^*_N(\mu),\mb{F}^*_N(\mu))$ satisfy the systems \cref{24-4-19ct3*}--\cref{24-4-19ct7*} and \cref{24-4-19ct3**}--\cref{24-4-19ct7**}  at a given $\mu\in \mathcal{P}$, respectively. Then the estimates
\begin{align*}
\underline{\delta}
&\le 
\|\mb{u}^*_h(\mu) - \mathbf{u}^*_N(\mu)\|_{\mathbf{L}^2(\Omega)} 
+ \|\mb{E}^*_h(\mu) - \mb{E}^*_N(\mu)\|_{\mathbf{H}(\mathbf{curl};\Omega)} \\
&~\quad + \|\mb{F}^*_h(\mu) - \mb{F}^*_N(\mu)\|_{\mathbf{H}(\mathbf{curl};\Omega)} 
\le \overline{\delta}
\end{align*}
are satisfied, where 
\begin{align*}
\overline{\delta} 
&= \overline{\delta}_{\mb{E}} \|R_{\mb{E}}(\cdot;\mu)\|_{{_{\mathbf{H}(\mathbf{curl};\Omega)}}^*} + \overline{\delta}_{\mb{F}} \|R_{\mb{F}}(\cdot;\mu)\|_{{_{\mathbf{H}(\mathbf{curl};\Omega)}}^*}\\
\overline{\delta}_{\mb{E}} 
&= C^\Omega_{\overline{\sigma}} \left( \alpha^{-1/2}\underline{\epsilon}^{-1}\overline{\epsilon}  + (1 + C^\Omega_{\overline{\sigma}} \overline{\epsilon}) \left(C^\Omega_{\overline{\sigma}}\alpha^{-1/2} \underline{\epsilon}^{-1} \overline{\epsilon}^2 + 1\right) \right) \\
\overline{\delta}_{\mb{F}} 
&= C^\Omega_{\overline{\sigma}} \big( 1 + \alpha^{-1}\underline{\epsilon}^{-1} \overline{\epsilon} + \left(1+ C^\Omega_{\overline{\sigma}} \overline{\epsilon} \right) C^\Omega_{\overline{\sigma}} \alpha^{-1}\underline{\epsilon}^{-1} \overline{\epsilon}^2 \big) 
\end{align*}
and
\begin{align*}
\underline{\delta} 
&= \underline{\delta}_{\mb{E}}\|R_{\mb{E}}(\cdot;\mu)\|_{{\mathbf{H}(\mathbf{curl};\Omega)}^*} + \underline{\delta}_{\mb{F}}\|R_{\mb{F}}(\cdot;\mu)\|_{{\mathbf{H}(\mathbf{curl};\Omega)}^*}\\
\underline{\delta}_{\mb{E}} 
&=  \underline{\sigma}\max(2, C^\Omega_{\overline{\sigma}} \overline{\epsilon})^{-1}\\
\underline{\delta}_{\mb{F}} 
&=  \underline{\sigma}\max(1, 2C^\Omega_{\overline{\sigma}} \overline{\epsilon})^{-1}.
\end{align*}
\end{theorem}

\begin{proof}
First we establish the upper bound. By variational inequalities \cref{24-4-19ct7*} and \cref{24-4-19ct7**}, we have
$$\left(\epsilon(\mathbf{u}^*_N - \mathbf{u}^*_h), \mathbf{u}_d  - \frac{1}{\alpha}\mathbf{F}^*_h - \mathbf{u}^*_h\right)_{\mathbf{L}^2(\Omega)} \le 0,~ \left(\epsilon(\mathbf{u}^*_h - \mathbf{u}^*_N), \mathbf{u}_d  - \frac{1}{\alpha}\mathbf{F}^*_N - \mathbf{u}^*_N\right)_{\mathbf{L}^2(\Omega)} \le 0$$
which implies that
\begin{align}
\alpha\underline{\epsilon}\|\mb{u}^*_h - \mathbf{u}^*_N\|^2_{\mathbf{L}^2(\Omega)} 
&\le \left(\epsilon(\mathbf{u}^*_N - \mathbf{u}^*_h), \mb{F}^*_h- \mathbf{F}^*_N\right)_{\mathbf{L}^2(\Omega)} \notag\\
&= \left(\epsilon(\mathbf{u}^*_N - \mathbf{u}^*_h), \mb{F}^*_h- \widehat{\mb{F}}_h\right)_{\mathbf{L}^2(\Omega)} + \left(\epsilon(\mathbf{u}^*_N - \mathbf{u}^*_h),\widehat{\mb{F}}_h - \mb{F}^*_N\right)_{\mathbf{L}^2(\Omega)}.
\label{2-5-19ct3}
\end{align}
Now we from \cref{29-4-19ct1}, \cref{24-4-19ct3*}, \cref{24-4-19ct5*} and \cref{29-4-19ct2} get
\begin{align}
&\left(\epsilon(\mathbf{u}^*_N - \mathbf{u}^*_h), \mb{F}^*_h- \widehat{\mb{F}}_h\right)_{\mathbf{L}^2(\Omega)} \notag\\
&~\quad = \left(\sigma^{-1}\nabla\times (\widehat{\mb{E}}_h - \mb{E}^*_h), \nabla\times (\mb{F}^*_h- \widehat{\mb{F}}_h)\right)_{\mathbf{L}^2(\Omega)}\notag\\
&~\quad = \left(\epsilon (\mathbf{E}^*_h - \mathbf{E}^*_N), \widehat{\mb{E}}_h - \mb{E}^*_h\right)_{\mathbf{L}^2(D)} \notag\\
&~\quad = \left(\epsilon(\mathbf{E}^*_h - \mathbf{E}^*_N), \widehat{\mb{E}}_h - \mb{E}^*_N + \mb{E}^*_N - \mb{E}^*_h\right)_{\mathbf{L}^2(D)} \notag\\
&~\quad \le \left(\epsilon(\mathbf{E}^*_h - \mathbf{E}^*_N), \widehat{\mb{E}}_h - \mb{E}^*_N \right)_{\mathbf{L}^2(D)}
- \underline{\epsilon}\left\|\mathbf{E}^*_h - \mathbf{E}^*_N\right\|^2_{\mathbf{L}^2(D)} \notag\\
&~\quad \le 2^{-1} \underline{\epsilon}^{-1}\overline{\epsilon}^2\| \widehat{\mb{E}}_h - \mb{E}^*_N \|^2_{\mathbf{H}(\mathbf{curl};\Omega)}
\label{2-5-19ct4}
\end{align}
and 
\begin{align} 
&\left(\epsilon(\mathbf{u}^*_N - \mathbf{u}^*_h),\widehat{\mb{F}}_h - \mb{F}^*_N\right)_{\mathbf{L}^2(\Omega)}\notag\\
&~\quad  \le 2^{-1}\alpha \underline{\epsilon}\|\mathbf{u}^*_N - \mb{u}^*_h\|^2_{\mathbf{L}^2(\Omega)} 
+ 2^{-1}\alpha^{-1} \underline{\epsilon}^{-1}\overline{\epsilon}^2\| \widehat{\mb{F}}_h - \mb{F}^*_N \|^2_{\mathbf{H}(\mathbf{curl};\Omega)}.
\label{2-5-19ct5}
\end{align}
We thus have from \cref{2-5-19ct3}--\cref{2-5-19ct5} that
\begin{align}
&\|\mb{u}^*_h - \mathbf{u}^*_N\|_{\mathbf{L}^2(\Omega)} \notag\\
&~\quad \le \alpha^{-1/2} \underline{\epsilon}^{-1} \overline{\epsilon} \| \widehat{\mb{E}}_h - \mb{E}^*_N \|_{\mathbf{H}(\mathbf{curl};\Omega)} + \alpha^{-1}\underline{\epsilon}^{-1} \overline{\epsilon}\| \widehat{\mb{F}}_h - \mb{F}^*_N \|_{\mathbf{H}(\mathbf{curl};\Omega)}.
\label{2-5-19ct6}
\end{align}
An application of \cref{30-4-19} and \cref{2-5-19ct6} yield
\begin{align}
\|\mb{E}^*_h - \mb{E}^*_N\|_{\mathbf{H}(\mathbf{curl};\Omega)} 
&\le \|\mb{E}^*_h - \widehat{\mb{E}}_h\|_{\mathbf{H}(\mathbf{curl};\Omega)} + \|\widehat{\mb{E}}_h - \mb{E}^*_N\|_{\mathbf{H}(\mathbf{curl};\Omega)} \notag\\
&\le C^\Omega_{\overline{\sigma}} \overline{\epsilon}\|\mb{u}^*_h - \mathbf{u}^*_N\|_{\mathbf{L}^2(\Omega)} + \|\widehat{\mb{E}}_h - \mb{E}^*_N\|_{\mathbf{H}(\mathbf{curl};\Omega)} \notag\\ 
&\le \left(C^\Omega_{\overline{\sigma}}\alpha^{-1/2} \underline{\epsilon}^{-1} \overline{\epsilon}^2 + 1\right) \| \widehat{\mb{E}}_h - \mb{E}^*_N \|_{\mathbf{H}(\mathbf{curl};\Omega)} \notag\\
&~\quad + C^\Omega_{\overline{\sigma}} \alpha^{-1} \underline{\epsilon}^{-1} \overline{\epsilon}^2\| \widehat{\mb{F}}_h - \mb{F}^*_N \|_{\mathbf{H}(\mathbf{curl};\Omega)}
\label{2-5-19ct7}
\end{align}
and
\begin{align}
\|\mb{F}^*_h - \mb{F}^*_N\|_{\mathbf{H}(\mathbf{curl};\Omega)}  
&\le \|\mb{F}^*_h - \widehat{\mb{F}}_h\|_{\mathbf{H}(\mathbf{curl};\Omega)} + \|\widehat{\mb{F}}_h - \mb{F}^*_N\|_{\mathbf{H}(\mathbf{curl};\Omega)} \notag\\
%&~\quad \le C^\Omega_{\overline{\sigma}} \overline{\epsilon}\|\mb{E}^*_h - \mathbf{E}^*_N \|_{\mathbf{L}^2(\Omega)} + \|\widehat{\mb{F}}_h - \mb{F}^*_N\|_{\mathbf{H}(\mathbf{curl};\Omega)} \notag\\
&~\quad \le C^\Omega_{\overline{\sigma}} \overline{\epsilon} \left(C^\Omega_{\overline{\sigma}}\alpha^{-1/2} \underline{\epsilon}^{-1} \overline{\epsilon}^2 + 1\right) \| \widehat{\mb{E}}_h - \mb{E}^*_N \|_{\mathbf{H}(\mathbf{curl};\Omega)} \notag\\
&~\quad\quad + \left( \left( C^\Omega_{\overline{\sigma}}\right)^2 \alpha^{-1} \underline{\epsilon}^{-1} \overline{\epsilon}^3 +1 \right)\| \widehat{\mb{F}}_h - \mb{F}^*_N \|_{\mathbf{H}(\mathbf{curl};\Omega)}.
\label{2-5-19ct8}
\end{align}
Combining the inequalities \cref{2-5-19ct6}--\cref{2-5-19ct8} with \cref{2-5-19}, we therefore arrive at the upper bound. Next, we will derive for the lower bound. We from \Cref{30-4-19} and \Cref{2-5-19} get that
\begin{align*}
\underline{\sigma} \|R_{\mb{E}}(\cdot;\mu)\|_{{_{\mathbf{H}(\mathbf{curl};\Omega)}}^*}
& \le \|\widehat{\mb{E}}_h - \mb{E}^*_h\|_{\mathbf{H}(\mathbf{curl};\Omega)} + \|\mb{E}^*_h - \mb{E}^*_N\|_{\mathbf{H}(\mathbf{curl};\Omega)} \notag\\
& \le C^\Omega_{\overline{\sigma}} \overline{\epsilon}\|\mb{u}^*_h - \mathbf{u}^*_N\|_{\mathbf{L}^2(\Omega)} + \|\mb{E}^*_h - \mb{E}^*_N\|_{\mathbf{H}(\mathbf{curl};\Omega)} \notag\\
& \le \max(2, C^\Omega_{\overline{\sigma}} \overline{\epsilon}) \left(\|\mb{u}^*_h - \mathbf{u}^*_N\|_{\mathbf{L}^2(\Omega)} + 2^{-1}\|\mb{E}^*_h - \mb{E}^*_N\|_{\mathbf{H}(\mathbf{curl};\Omega)} \right)
\end{align*}
and
\begin{align*} 
\underline{\sigma} \|R_{\mb{F}}(\cdot;\mu)\|_{{_{\mathbf{H}(\mathbf{curl};\Omega)}}^*} 
& \le \|\widehat{\mb{F}}_h - \mb{F}^*_h\|_{\mathbf{H}(\mathbf{curl};\Omega)} + \|\mb{F}^*_h - \mb{F}^*_N\|_{\mathbf{H}(\mathbf{curl};\Omega)} \notag\\ 
& \le C^\Omega_{\overline{\sigma}} \overline{\epsilon}\|\mb{E}^*_h - \mathbf{E}^*_N \|_{\mathbf{H}(\mathbf{curl};\Omega)} + \|\mb{F}^*_h - \mb{F}^*_N\|_{\mathbf{H}(\mathbf{curl};\Omega)} \notag\\
& \le \max(1, 2C^\Omega_{\overline{\sigma}} \overline{\epsilon}) \left(\frac{1}{2}\|\mb{E}^*_h - \mathbf{E}^*_N \|_{\mathbf{H}(\mathbf{curl};\Omega)} + \|\mb{F}^*_h - \mb{F}^*_N\|_{\mathbf{H}(\mathbf{curl};\Omega)}\right)
\end{align*}
which completes the proof.
\end{proof}

We aim towards an error bound for the cost functional.

\begin{theorem}\label{2-5-19-DL*}
Let $(\mb{u}^*_h(\mu),\mb{E}^*_h(\mu),\mb{F}^*_h(\mu))$ and $(\mb{u}^*_N(\mu),\mb{E}^*_N(\mu),\mb{F}^*_N(\mu))$ satisfy the systems \cref{24-4-19ct3*}--\cref{24-4-19ct7*} and \cref{24-4-19ct3**}--\cref{24-4-19ct7**}  at a given $\mu\in \mathcal{P}$, respectively. Then, 
\begin{align*}
\left| J(\mb{u}^*_h,\mb{E}^*_h;\mu) - J(\mb{u}^*_N,\mb{E}^*_N;\mu) \right| \le \delta^J,
\end{align*}
with
\begin{align*}
\delta^J 
&= \delta^J_{\mb{E}} \|R_{\mb{E}}(\cdot;\mu)\|_{{_{\mathbf{H}(\mathbf{curl};\Omega)}}^*} + \delta^J_{\mb{F}} \|R_{\mb{F}}(\cdot;\mu)\|_{{_{\mathbf{H}(\mathbf{curl};\Omega)}}^*}\\
\delta^J_{\mb{E}} 
&= C^\Omega_{\overline{\sigma}} \overline{\epsilon} \left( (C_\mathbf{E} + e_d) \left(C^\Omega_{\overline{\sigma}}\alpha^{-1/2} \underline{\epsilon}^{-1} \overline{\epsilon}^2 + 1\right) + \alpha^{1/2} \underline{\epsilon}^{-1} \overline{\epsilon}\left( \|\overline{\mb{u}}\|_{\mathbb{R}^3} |\Omega|^{1/2} + u_d \right)\right)\\
\delta^J_{\mb{F}}
&= C^\Omega_{\overline{\sigma}} \overline{\epsilon}^2 \left( C^\Omega_{\overline{\sigma}} \alpha^{-1} \underline{\epsilon}^{-1} \overline{\epsilon} (C_\mathbf{E} + e_d) + \underline{\epsilon}^{-1}  \left( \|\overline{\mb{u}}\|_{\mathbb{R}^3} |\Omega|^{1/2} + u_d \right)\right).
\end{align*}
\end{theorem}

\begin{proof}
We get that
\begin{align}\label{10-5-19ct3}
&2\left| J(\mb{u}^*_h,\mb{E}^*_h;\mu) - J(\mb{u}^*_N,\mb{E}^*_N;\mu) \right| \notag\\
&~\quad = \left( \epsilon (\mb{E}^*_h - \mb{E}^*_N), \mb{E}^*_h + \mb{E}^*_N- 2\mb{E}_d\right)_{\mb{L}^2(D)} \notag\\
&~\quad\quad + \alpha \left( \epsilon (\mb{u}^*_h - \mb{u}^*_N), \mb{u}^*_h + \mb{u}^*_N- 2\mb{u}_d\right)_{\mb{L}^2(\Omega)}\notag\\
&~\quad \le \overline{\epsilon}\left( \|\mb{E}^*_h\|_{\mb{L}^2(\Omega)}+ \|\mb{E}^*_N\|_{\mb{L}^2(\Omega)} + 2\|\mb{E}_d\|_{\mb{L}^2(D)}\right) \|\mb{E}^*_h - \mb{E}^*_N\|_{\mb{L}^2(\Omega)}\notag\\
&~\quad\quad + \overline{\epsilon} \alpha \left( \|\mb{u}^*_h\|_{\mb{L}^2(\Omega)}+ \|\mb{u}^*_N\|_{\mb{L}^2(\Omega)} + 2\|\mb{u}_d\|_{\mb{L}^2(\Omega)}\right) \|\mb{u}^*_h - \mb{u}^*_N\|_{\mb{L}^2(\Omega)}\notag\\
&~\quad \le 2\overline{\epsilon} \left( C_\mathbf{E} +e_d\right) \|\mb{E}^*_h - \mb{E}^*_N\|_{\mb{L}^2(\Omega)}\notag\\
&~\quad\quad + 2\overline{\epsilon} \alpha \left( \|\overline{\mb{u}}\|_{\mathbb{R}^3} |\Omega|^{1/2} + u_d \right) \|\mb{u}^*_h - \mb{u}^*_N\|_{\mb{L}^2(\Omega)}.
\end{align}
By \cref{2-5-19ct7}, \cref{2-5-19ct6} and \cref{2-5-19}, we have
\begin{align}\label{10-5-19ct1}
\|\mb{E}^*_h - \mb{E}^*_N\|_{\mathbf{H}(\mathbf{curl};\Omega)} 
&\le C^\Omega_{\overline{\sigma}} \left(C^\Omega_{\overline{\sigma}}\alpha^{-1/2} \underline{\epsilon}^{-1} \overline{\epsilon}^2 + 1\right) \|R_{\mb{E}}\|_{{_{\mathbf{H}(\mathbf{curl};\Omega)}}^*} \notag\\
&~\quad + \left( C^\Omega_{\overline{\sigma}}\right)^2 \alpha^{-1} \underline{\epsilon}^{-1} \overline{\epsilon}^2 \|R_{\mb{F}}\|_{{_{\mathbf{H}(\mathbf{curl};\Omega)}}^*}
\end{align}
and 
\begin{align}\label{10-5-19ct2}
\|\mb{u}^*_h - \mathbf{u}^*_N\|_{\mathbf{L}^2(\Omega)} 
&\le C^\Omega_{\overline{\sigma}}\alpha^{-1/2} \underline{\epsilon}^{-1} \overline{\epsilon} \|R_{\mb{E}}\|_{{_{\mathbf{H}(\mathbf{curl};\Omega)}}^*} \notag\\
&~\quad + C^\Omega_{\overline{\sigma}} \alpha^{-1} \underline{\epsilon}^{-1} \overline{\epsilon} \|R_{\mb{F}}\|_{{_{\mathbf{H}(\mathbf{curl};\Omega)}}^*}.
\end{align}
The desired inequality follows directly from the estimates \cref{10-5-19ct3}--\cref{10-5-19ct2}, which finishes the proof.
\end{proof}

We derive a posteriori error estimators.

\begin{theorem}\label{3-5-19ct6}
The absolute a posteriori error estimator
\begin{align*}
\|\mb{u}^*_h(\mu) - \mathbf{u}^*_N(\mu)\|_{\mathbf{L}^2(\Omega)} 
\le \Delta^{\mbox{ab}}_N(\mu)
\end{align*}
is established, where
\begin{align*}
\Delta^{\mbox{ab}}_N(\mu) 
:=  C^\Omega_{\overline{\sigma}}\alpha^{-1/2} \underline{\epsilon}^{-1} \overline{\epsilon} \|R_{\mb{E}}(\cdot;\mu)\|_{{_{\mathbf{H}(\mathbf{curl};\Omega)}}^*}  + C^\Omega_{\overline{\sigma}} \alpha^{-1} \underline{\epsilon}^{-1} \overline{\epsilon} \|R_{\mb{F}}(\cdot;\mu)\|_{{_{\mathbf{H}(\mathbf{curl};\Omega)}}^*}.
\end{align*}
Furthermore, in case $\frac{2\Delta^{\mbox{ab}}_N}{\|\mathbf{u}^*_N\|_{\mathbf{L}^2(\Omega)}} \le 1$ we have the relative a posteriori error estimator
\begin{align}\label{3-5-19ct5}
\dfrac{\|\mb{u}^*_h(\mu) - \mathbf{u}^*_N(\mu)\|_{\mathbf{L}^2(\Omega)}}{\|\mathbf{u}^*_h(\mu)\|_{\mathbf{L}^2(\Omega)}}
\le \Delta^{\mbox{re}}_N(\mu):= \dfrac{2\Delta^{\mbox{ab}}_N(\mu)}{\|\mathbf{u}^*_N(\mu)\|_{\mathbf{L}^2(\Omega)}}.
\end{align}
\end{theorem}

\begin{proof}
By \cref{10-5-19ct2}, it remains to show \cref{3-5-19ct5} only. We get
$$\|\mathbf{u}^*_N\|_{\mathbf{L}^2(\Omega)} - \|\mathbf{u}^*_h\|_{\mathbf{L}^2(\Omega)}
\le \|\mb{u}^*_h - \mathbf{u}^*_N\|_{\mathbf{L}^2(\Omega)} 
\le \Delta^{\mbox{ab}}_N
\le 2^{-1}\|\mathbf{u}^*_N\|_{\mathbf{L}^2(\Omega)}$$
and \cref{3-5-19ct5} is derived, which finishes the proof.
\end{proof}

%\section*{Acknowledgments}
%We would like to acknowledge the assistance of volunteers in putting
%together this example manuscript and supplement.

\end{document}